\documentclass[12pt]{amsart}

\calclayout

\usepackage[utf8]{inputenc}
\usepackage{amsmath,amsthm,amssymb}
\usepackage{colonequals}
\usepackage{enumitem}
\usepackage{mathtools}

\usepackage[colorlinks=true,allcolors=blue,pdftitle={Clearing out}]{hyperref}
\usepackage[initials,nobysame,alphabetic,msc-links,backrefs]{amsrefs}

\newtheorem{theorem}{Theorem}[section]
\newtheorem{theoremA}{Theorem}

\newtheorem{lemma}{Lemma}[section]
\newtheorem{lemmaA}{Lemma}

\newtheorem{prop}{Proposition}[section]
\newtheorem{propA}{Proposition}

\newtheorem{remark}{Remark}[section]
\newtheorem{definition}{Definition}[section]

\numberwithin{equation}{section}

\DeclareMathOperator{\diver}{div}
\DeclareMathOperator{\curl}{curl}

\newcommand{\R}{\ensuremath{\mathbb{R}}}
\newcommand{\N}{\ensuremath{\mathbb{N}}}
\newcommand{\C}{\ensuremath{\mathbb{C}}}
\newcommand{\norm}[3][]{#1\left\|#2#1\right\|_{#3}}
\newcommand{\ip}[2]{\left\langle #1\hspace*{.5mm},#2\right\rangle}
\newcommand{\dist}{\mathrm{dist}}

\newcommand{\ep}{\varepsilon}
\newcommand{\ex}{\mathrm{ex}}
\newcommand{\loc}{\mathrm{loc}}
\newcommand{\lep}{{|\log\ep|}}

\renewcommand{\u}{{\mathsf{u}}}
\newcommand{\A}{{\mathsf{A}}}
\renewcommand{\H}{{\mathsf{H}}}

\newcommand{\pmin}{\ensuremath{\rho_\ep}}
\newcommand{\Bfield}{\ensuremath{B_\ep^0}}

\newcommand{\fenthreedOmega}{\ensuremath{F_{\ep,\Omega}}}
\newcommand{\fenweightthreed}{\ensuremath{F_{\ep,\pmin}}}
\newcommand{\glthreed}{\ensuremath{GL_\ep}}

\newcommand{\phase}{\ensuremath{e^{ih_\ex\phi_\ep^0}}}
\newcommand{\meisconf}{\ensuremath{\pmin\phase,h_\ex A_\ep^0}}

\newcommand{\ACouGauge}{\ensuremath{\A^C}}
\newcommand{\AUniGauge}{\ensuremath{\widetilde{\A}}}
\newcommand{\uCouGauge}{\ensuremath{\u^C}}

\newcommand{\varphiCouGauge}{\ensuremath{\varphi^C}}

\newcommand{\AfluctWeight}{\ensuremath{\widetilde{\A}_{\mathrm{fl}}}}
\newcommand{\Afluctglobfree}{\ensuremath{\widetilde{\A}^{\mathrm{fl}}}}

\newcommand{\critfield}{\ensuremath{H^{\varepsilon}_{c_1}}}
\newcommand{\pzero}{{p_0}}
\newcommand{\minsp}{H^1(\Omega,\C)\times[A_{\ex}+H_{\curl}]}
\newcommand{\McurlBfield}{\ensuremath{\widetilde{M}_s}}

\title[Clearing-out for the pinned three-dimensional Ginzburg--Landau model]{Clearing-out for the pinned three-dimensional Ginzburg--Landau model}

\author[Mat\'ias D\'iaz-Vera]{Mat\'ias D\'iaz-Vera}
\address[Mat\'ias D\'iaz-Vera]{School of Mathematics, Georgia Institute of Technology, Atlanta, GA 30332-0160, United States of America}
\email{mdiaz88@gatech.edu}

\author[Carlos Rom\'an]{Carlos Rom\'an}
\address[Carlos Rom\'an]{Facultad de Matem\'aticas \& Instituto de Ingenier\'ia Matem\'atica y Computacional, Pontificia Universidad Cat\'olica de Chile, Vicu\~na Mackenna 4860, 7820436 Macul, Santiago, Chile}
\email{carlos.roman@uc.cl}

\author[Kai Yamamoto]{Kai Yamamoto}
\address[Kai Yamamoto]{Instituto de Ingenier\'ia Matem\'atica y Computacional, Pontificia Universidad Cat\'olica de Chile, Vicu\~na Mackenna 4860, 7820436 Macul, Santiago, Chile}
\email{kai.yamamoto@uc.cl}

\date{November 27, 2026}



\thanks{Funding information: ANID FONDECYT 1231593.}

\begin{document}
	
\begin{abstract}
We continue our study of global minimizers of the three-dimensional magnetic Ginzburg--Landau energy with a spatially dependent pinning term. Our main result is a clearing-out theorem that upgrades weak vortexlessness, previously obtained through energy and vorticity estimates, to a strong pointwise characterization of the Meissner regime. More precisely, we prove that, below the first critical field, the modulus of the order parameter remains uniformly close to the pinning profile. Unlike in the homogeneous model, the presence of a spatially dependent pinning profile generates additional terms in the Pohozaev identity, making the derivation of the pseudomonotonicity formula underlying the clearing-out argument substantially more delicate.
\end{abstract}
\maketitle 

\section{Introduction}
The magnetic Ginzburg--Landau model is one of the fundamental mathematical frameworks for the study of superconductivity. Introduced by Ginzburg and Landau \cite{GinLan}, it provides a phenomenological description of superconducting materials and has strong connections with other variational models, including the Gross--Pitaevskii equation for superfluids \cite{bethuel-orlandi-smets-grosspitaevski} and models of rotating Bose--Einstein condensates \cite{aftalionbook}. 

\medskip
At the heart of the theory lies the \emph{Ginzburg--Landau energy functional}. In this work, we consider its inhomogeneous version, in which a spatially dependent pinning term models variations in the superconducting material:
\begin{equation}\label{def: gl energy no dec}
	\glthreed(\u,\A)\colonequals\frac12\int_\Omega\left(|\nabla_\A\u|^2+\frac1{2\ep^2}(a_\ep-|\u|^2)^2\right)+\frac12\int_{\R^3}|\H-\H_\ex|^2.
\end{equation}
Minimizing configurations $(\u,\A)$ of \eqref{def: gl energy no dec} describe the equilibrium state of a superconductor exposed to the external magnetic field $H_\ex$. Here:
\begin{itemize}
	\item $\Omega\subset\R^3$ is a bounded simply connected domain with boundary of class $C^2$, representing the superconducting sample.
	\item The configuration $(\u,\A)$ consists of a complex-valued function $\u$ defined in $\Omega$, whose modulus squared describes the local density of superconducting electrons, and a vector field $\A$ defined in $\R^3$, whose curl gives the induced magnetic field $\H$, that is, $\H=\curl\A$.
	\item $\nabla_\A\colonequals\nabla-i\A$ is the covariant derivative, so that $\nabla_\A\u$ is a complex-valued vector field.
	\item The pinning term $a_\ep$ is a measurable real-valued function taking values in $[b,1]$, where $b\in(0,1)$ is independent of $\ep$. It describes the local equilibrium density of superconducting electrons. A nonconstant $a_\ep$ models spatial inhomogeneities in the superconducting material and gives rise to vortex pinning.
	\item The vector field $H_\ex:\R^3\to\R^3$ represents the applied magnetic field. We assume that $H_\ex=h_\ex H_{0,\ex}$, where $H_{0,\ex}$ is a fixed vector field describing its spatial profile, normalized so that $\|H_{0,\ex}\|_{L^\infty(\R^3,\R^3)}=1$, and $h_\ex\geq0$ is an adjustable scalar parameter describing its intensity.
	\item The parameter $\ep>0$ is the inverse of the \emph{Ginzburg--Landau parameter}, usually denoted by $\kappa$, a dimensionless parameter depending only on the material. 
\end{itemize}

The magnetic Ginzburg--Landau model is a $\mathbb U(1)$-gauge theory, meaning that the Ginzburg--Landau energy is invariant under the gauge transformations
$$
	\u\to\u e^{i\phi},\qquad \A\to\A+\nabla\phi,
$$
where $\phi$ is an arbitrary sufficiently regular function. Physically relevant quantities remain invariant under these transformations, including the \emph{free energy} 
$$
	F_\ep(\u,\A)\colonequals\frac12\int_\Omega\left(|\nabla_\A\u|^2+\frac1{2\ep^2}(a_\ep-|\u|^2)^2\right)+\frac12\int_{\R^3}|\curl\A|^2,
$$
and the \emph{vorticity measure}, which is the gauge-invariant version of twice the Jacobian,
$$
	\mu(\u,\A)\colonequals\curl\left(\ip{i\u}{\nabla_\A\u}+\A\right).
$$
Here, $\ip{\cdot}{\cdot}$ denotes the real inner product in $\C$, namely, $\ip{z}{w}\colonequals\Re(\overline zw)$, applied componentwise when taking the inner product between a complex scalar and a complex vector.

\medskip
Superconductors are classified as type II when $\kappa>1/\sqrt{2}$, or equivalently, when $\ep<\sqrt{2}$. They are characterized by the existence of a mixed phase in which superconductivity coexists with \emph{vortices}, codimension-two singularities where superconductivity is locally suppressed. In the Ginzburg--Landau model, vortices are associated with the zero set of the order parameter and are characterized by the nonzero integer degree of its phase around them. In this work, we focus on the regime $\ep\ll1$, corresponding to \emph{extreme type-II superconductors}.

For sufficiently weak applied magnetic fields, the superconducting sample remains free of vortices. This vortexless regime is called the \emph{Meissner state}. As the intensity $h_\ex$ of the applied field increases, the formation of vortices eventually becomes energetically favorable. The field intensity separating these two regimes is called the \emph{first critical field} and is denoted by $H_{c_1}$.

\medskip

The inhomogeneous version of the Ginzburg--Landau model dates back to the work of Likharev \cite{Likha} and provides a more physically accurate description of superconducting materials. The relevance of inhomogeneities lies in their ability to restrict the motion of vortices by pinning them to preferred regions of the material, thereby reducing the energy dissipation caused by their motion. This explains why $a_\ep$ is referred to as a pinning term.

The pinning phenomenon is of significant interest in the physics community. Recent work such as \cite{Kwok2016} discusses the effect of different types of defects, which give rise to complex pinning landscapes influencing vortex configurations and enhancing the maximal current that can be carried without losing superconductivity. In a more microscopic context, experiments on FeSe-based superconductors exhibit local defects acting as pinning regions \cite{phys-x-2024}. This motivates the study of models incorporating material defects in superconducting samples. We also refer to \cite{Tin} for further physics-oriented discussions on vortex pinning.

\medskip
The mathematical analysis of the Ginzburg--Landau model began with the homogeneous case, in which $a_\ep\equiv1$. Its theory plays a major role in our work, although our setting allows for a broad class of spatially varying pinning terms. The seminal work of Bethuel, Brezis, and Hélein \cite{BetBreHel} laid the foundations for the asymptotic analysis of vortices and initiated an extensive mathematical literature. We refer to \cite{SanSerBook} and the references therein for a comprehensive account of the two-dimensional theory.

In three dimensions, vortices have a filamentary structure, which introduces substantial geometric difficulties compared with the point-vortex structure in two dimensions. Relevant works such as \cites{Riv,LinRiviere,sandier-surface-ball-construction,AlaBroMon,BalJerOrlSon1,BalJerOrlSon2} provide qualitative results on the limiting vortex structure of minimizers in different regimes. The study of minimizers near the first critical field, however, requires more refined tools that are quantitative in $\ep$. Such tools were developed in \cite{roman-vortex-construction}; we refer to \cites{Rom-CMP,roman-sandier-serfaty,RSS2} for further developments and applications.

\medskip
The literature on two-dimensional models including a pinning term is extensive. Some early works \cites{chapman-richardson,chapman-du-gunzburguer-variable-thickness,lin-du} study the effect of locally varying thickness on the location of vortices. More mathematically oriented analyses of specific pinning models, including smooth or periodic structures for $a_\ep$, can be found in \cites{ass,kachmar,dos-santos-2019,DosRodSan}. A characterization of the first critical field for a broad class of pinning terms $a_\ep$ is provided in \cite{diaz-roman-2d}.

Compared with the two-dimensional problem, the literature on the three-dimensional pinning model is considerably less developed. The work \cite{DosSantos3D} describes the limiting vortex set in the simpler setting without magnetic components, where the vortex structure is induced by a Dirichlet boundary condition. A section of \cite{DosRodSan} considers \eqref{def: gl energy no dec} when $a_\ep$ has a stationary ergodic structure at scale $\ep$, leading, as $\ep\to0$, to the homogeneous limiting model studied in \cite{BalJerOrlSon1}.

\medskip
In \cite{diaz-roman-3d}, the lower bound for the leading order of the first critical field in the three-dimensional inhomogeneous problem was first addressed, characterizing the Meissner state through a \emph{weak} notion of vortexlessness, in the sense that both the vorticity measure and the excess free energy relative to the Meissner state vanish as $\ep\to0$. A key ingredient in this analysis was the development of an inhomogeneous version of the $\ep$-level estimates from \cite{roman-vortex-construction}. A complementary upper bound for $H_{c_1}$ was obtained in \cite{RomanPinningUpper25}. More recently, the metastability and uniqueness of the Meissner state for applied fields much larger than $H_{c_1}$, reaching regimes of order $\ep^{-\alpha}$ for $\alpha<1/4$, were studied in \cite{rom2026metastab}.

\medskip
A natural question is whether this weak form of vortexlessness can be upgraded to a pointwise characterization of the Meissner state. More precisely, one would like to show that the modulus of the order parameter remains uniformly close to the pinning profile and therefore stays uniformly away from zero. In the three-dimensional setting, obtaining such a result requires a clearing-out theorem adapted to the weighted Ginzburg--Landau energy. Developing this tool and deriving the resulting strong characterization of the Meissner state below the first critical field are the main objectives of the present work.

\subsection{Decoupling of the pinning term and weighted Ginzburg--Landau reformulation}

A fundamental tool in the analysis of inhomogeneous Ginzburg--Landau models is the decoupling procedure introduced by Lassoued and Mironescu \cite{lassoued-mironescu}, which isolates the energetic contribution of the pinning term and reduces \eqref{def: gl energy no dec} to a weighted Ginzburg--Landau functional.

Specifically, let $\pmin$ denote the unique positive minimizer of the scalar energy
$$
E_\ep(v)\colonequals\frac12\int_\Omega\left(|\nabla v|^2+\frac1{2\ep^2}(a_\ep-|v|^2)^2\right),
$$
which solves the Euler--Lagrange equation
\begin{equation}\label{pde rho}
	\left\{
	\begin{array}{rcll}
		-\Delta\pmin&=&\dfrac{\pmin}{\ep^2}(a_\ep-\pmin^2)&\mathrm{in}\ \Omega\\[2mm]
		\dfrac{\partial\pmin}{\partial\nu}&=&0&\mathrm{on}\ \partial\Omega.
	\end{array}
	\right.
\end{equation}
We recall from \cite[Proposition 2.3]{RomanPinningUpper25} that
$$
b\leq \pmin^2\leq 1.
$$
Writing any admissible configuration $(\u,\A)$ in the form $\u=\pmin u$ and integrating by parts using \eqref{pde rho}, one obtains the decoupling identity
$$
	\glthreed(\pmin u,\A)=E_\ep(\pmin)+\frac12\int_\Omega\left(\pmin^2|\nabla_\A u|^2+\frac{\pmin^4}{2\ep^2}(1-|u|^2)^2\right)+\frac12\int_{\R^3}|\H-H_\ex|^2.
$$
Since $E_\ep(\pmin)$ depends only on $a_\ep$ and $\Omega$, the analysis of \eqref{def: gl energy no dec} reduces to studying the weighted energy
\begin{equation}\label{eq:weighted gl energy}
	\frac12\int_\Omega\left(\pmin^2|\nabla_\A u|^2+\frac{\pmin^4}{2\ep^2}(1-|u|^2)^2\right)+\frac12\int_{\R^3}|\H-H_\ex|^2.
\end{equation}

The profile $\pmin^2$ can be viewed as a regularized version of $a_\ep$, as $E_\ep$ penalizes both strong oscillations and deviations from $a_\ep$. In two dimensions, the pinning mechanism is particularly transparent in the parabolic Ginzburg--Landau model \cite[Remark 2.1]{duerinckx-serfaty}, where $-\nabla\log\pmin^2$ acts as an effective force field driving vortices toward local minima of $\pmin$.

\subsection{Clearing-out and strong vortexlessness}

One of the key tools in the characterization of minimizing configurations is the \emph{clearing-out property}, also referred to as $\eta$-ellipticity or $\eta$-compactness in earlier literature. Roughly speaking, it asserts that if the local free energy of a configuration in $B_r(x_0)\cap\Omega$, normalized by the radius $r$, is sufficiently small relative to $|\log\ep|$, then the configuration is locally vortexless in the strong sense that the modulus of its order parameter is uniformly bounded away from zero.

The first clearing-out result was established in two dimensions by Bethuel, Brezis, and Hélein \cite[Theorem III.3]{BetBreHel} through elementary estimates. Specifically, given the standard gradient bound of order $\ep^{-1}$ available in the relevant regime, any point where $|u|$ stays away from $1$ forces $|u|$ to remain away from $1$ on a ball of radius comparable to $\ep$. Since the area of such a ball cancels the $\ep^{-2}$ factor in the potential term, the potential energy is bounded below by a positive constant independent of $\ep$. Consequently, if the energy is sufficiently small, zeros are immediately excluded without directly invoking minimality or the Euler--Lagrange equations, which allows the argument to be extended directly to inhomogeneous two-dimensional models, as in \cite[Proposition 2.5]{diaz-roman-2d}.

In three dimensions, the first clearing-out results were proved by Lin and Rivière \cites{Riv,LinRiviere}, and were subsequently extended to non-minimizing solutions and up to the boundary in \cites{BetBreOrl,bethuel-orlandi-smets-grosspitaevski,chiron2005boundary}. In particular, the boundary clearing-out estimate of Chiron \cite[Theorem 3]{chiron2005boundary} was a key ingredient in \cite{Rom-CMP} to establish strong vortexlessness below $H_{c_1}$ in the homogeneous three-dimensional model. Because no comparable tool was available in the inhomogeneous setting, the analysis in \cite{diaz-roman-3d} was restricted to weak vortexlessness. Developing a clearing-out result adapted to the inhomogeneous model is therefore essential to match the precision of the homogeneous theory.

In contrast with the two-dimensional case, proving clearing-out in three dimensions is considerably more involved and relies essentially on the Euler--Lagrange equation for the order parameter $u$. Although the elementary scaling argument of \cite{BetBreHel} still applies in three dimensions, the volume of a ball of radius comparable to $\ep$ is of order $\ep^3$, so a point where $|u|$ stays away from $1$ only yields a lower bound of order $\ep$ on the potential energy; thus, ruling out zeros requires smallness of the normalized energy at scale $\ep$ (see Lemma~\ref{lemma:clearingout}). To bridge macroscopic scales and scale $\ep$, the standard strategy is to derive a monotonicity formula for a normalized local energy by means of a Pohozaev identity, obtained by testing the equation for $u$ with $(x-x_0)\cdot\nabla u$. In the magnetic setting, fixing a suitable gauge and expanding the covariant Laplacian allows one to rewrite the equation for $u$ in terms of the standard Laplacian: the first-order magnetic coupling is treated as a lower-order perturbation in the Pohozaev identity, while the quadratic magnetic potential term is absorbed into the potential well, as in \cite{chiron2005boundary}. This leads to a normalized local energy of the form
$$
	\frac1{2r}\int_{B_r(x_0)\cap\Omega}\left(|\nabla u|^2+\frac1{2\ep^2}(1-\ep^2|\A|^2-|u|^2)^2\right).
$$
The resulting monotonicity formula propagates energy bounds from macroscopic scales down to scale $\ep$, where the elementary scaling argument yields uniform smallness of $1-|u|$ and hence rules out zeros.

In our inhomogeneous setting, applying the same reformulation to the Euler--Lagrange equation for $u$ associated with \eqref{eq:weighted gl energy} leads us to consider the weighted normalized local energy
$$
	\frac1{2r}\int_{B_r(x_0)\cap\Omega}\left(\pmin^2|\nabla u|^2+\frac{\pmin^4}{2\ep^2}(b_\ep-|u|^2)^2\right),
$$
with the modified coefficient
$$
b_\ep\colonequals1-\ep^2\frac{|\A|^2}{\pmin^2}.
$$
The presence of the weights $\pmin^2$ and $\pmin^4$ introduces a significant difficulty: derivatives of $\pmin$ appear in the Pohozaev identity and generate additional error terms that cannot be controlled by the local energy under a mere Hölder bound on $\pmin$. Carrying out this argument therefore requires a quantitative Lipschitz estimate on the pinning profile.

Our clearing-out proof follows the general framework of \cite[Theorem 2]{BetBreOrl} and incorporates boundary techniques from \cite[Theorem 3]{chiron2005boundary}. We note that a clearing-out estimate with pinning was previously considered in \cite[Theorem 3.1]{LiuZhou}; however, that work is restricted to a smooth, $\ep$-independent pinning weight and to a reduced model in which the magnetic energy is confined to $\Omega$. In contrast, our analysis addresses the full three-dimensional magnetic energy in $\R^3$ and allows for general $\ep$-dependent pinning terms satisfying quantitative regularity bounds as $\ep\to0$. The main technical step consists in deriving a weighted Pohozaev identity, carefully controlling both interior and boundary error terms, and deducing a pseudomonotonicity formula that propagates energy bounds down to scale $\ep$.

\subsection{Main results}

To state our main results, we first specify the functional setting and introduce the weighted Meissner state. Since $\diver H_{0,\ex}=0$ in $\R^3$, in accordance with the nonexistence of magnetic monopoles in Maxwell's theory of electromagnetism, there exists a vector potential $A_{0,\ex}\in H^1_{\loc}(\R^3,\R^3)$ such that
$$
\curl A_{0,\ex}=H_{0,\ex}\quad\mathrm{in}\ \R^3.
$$
We define
$$
A_{\ex}\colonequals h_{\ex}A_{0,\ex}.
$$
Observe that we can choose $A_{0,\ex}$ in the Coulomb gauge in $\Omega$, namely,
$$
\left\{
\begin{array}{rcll}
	\diver A_{0,\ex}&=&0&\mathrm{in}\ \Omega\\
	A_{0,\ex}\cdot\nu&=&0&\mathrm{on}\ \partial\Omega,
\end{array}
\right.
$$
We minimize \eqref{def: gl energy no dec} over the admissible space $H^1(\Omega,\C)\times[A_\ex+H_{\curl}]$, where
$$
H_{\curl}\colonequals\left\{A\in H^1_{\loc}(\R^3,\R^3)\ \middle|\ \curl A\in L^2(\R^3,\R^3)\right\}.
$$
Let $\dot H^1(\R^3,\R^3)$ denote the completion of $C_c^\infty(\R^3,\R^3)$ with respect to the norm $\|\nabla(\,\cdot\,)\|_{L^2(\R^3,\R^3)}$, and set
$$
\dot H^1_{\diver=0}\colonequals\left\{A\in\dot H^1(\R^3,\R^3)\ \middle|\ \diver A=0\ \mathrm{in}\ \R^3\right\}.
$$
For every $A\in\dot H^1_{\diver=0}$, we have
\begin{equation}\label{space2}
	\|A\|_{\dot H^1_{\diver=0}}\colonequals\|\nabla A\|_{L^2(\R^3,\R^3)}=\|\curl A\|_{L^2(\R^3,\R^3)}.
\end{equation}

In view of the weighted reformulation \eqref{eq:weighted gl energy}, vortexless configurations are naturally approximated by unit-modulus order parameters $u=e^{i\phi}$. Rescaling by the intensity of the applied field, for any vector potential $A\in A_{0,\ex}+\dot H^1_{\diver=0}$ we choose the phase $\phi_A$ minimizing $\frac12\int_\Omega\pmin^2|A-\nabla\phi|^2$ among zero average functions, so that $A-\nabla\phi_A$ satisfies the weighted Coulomb gauge
\begin{equation}\label{eqphi0}
	\left\{
	\begin{array}{rcll}
		\diver\big(\pmin^2(A-\nabla\phi_A)\big)&=&0&\mathrm{in}\ \Omega\\
		(A-\nabla\phi_A)\cdot\nu&=&0&\mathrm{on}\ \partial\Omega.
	\end{array}
	\right.
\end{equation}
Minimizing the resulting reduced functional
\begin{equation}\label{reducedenergy}
	J_\ep(A)\colonequals\frac12\int_\Omega\pmin^2|A-\nabla\phi_A|^2+\frac12\int_{\R^3}|\curl(A-A_{0,\ex})|^2
\end{equation}
yields a unique minimizer $A_\ep^0$ satisfying
\begin{equation}\label{spaceforA0}
	A_\ep^0-A_{0,\ex}\in\dot H^1_{\diver=0}.
\end{equation}
Setting $\phi_\ep^0\colonequals\phi_{A_\ep^0}$, we refer to $(\meisconf)$ as the weighted Meissner state. Moreover, by \eqref{eqphi0}, there exists a unique vector field $\Bfield$ such that
\begin{equation}\label{defBfield}
	\pmin^2(A_\ep^0-\nabla\phi_\ep^0)=\curl\Bfield\quad\mathrm{in}\ \Omega,
\end{equation}
together with
\begin{equation}\label{conditionsB0}
	\left\{
	\begin{array}{rcll}
		\diver\Bfield&=&0&\mathrm{in}\ \Omega\\
		\Bfield\times\nu&=&0&\mathrm{on}\ \partial\Omega.
	\end{array}
	\right.
\end{equation}

\medskip
The leading order of the first critical field is determined by a weighted version of the \emph{isoflux problem}. Let $\mathcal{N}$ denote the space of normal $1$-currents supported in $\overline\Omega$ whose boundaries are supported in $\partial\Omega$, that is, $1$-currents of finite mass whose boundaries also have finite mass, and denote the mass of a current by $|\cdot|$. Identifying $\Bfield$ with a $1$-form, we denote by $\ip{\Gamma}{\Bfield}$ the duality pairing between $\Gamma\in\mathcal{N}$ and $\Bfield$, which coincides with the circulation of $\Bfield$ along $\Gamma$ whenever $\Gamma$ is a simple oriented Lipschitz curve. Moreover, the weighted current $\pmin^2\Gamma\in\mathcal{N}$ is defined by duality via $\ip{\pmin^2\Gamma}{\omega}\colonequals\ip{\Gamma}{\pmin^2\omega}$ for every $1$-form $\omega$ supported in $\Omega$.

We then define the weighted isoflux ratio
$$
	\mathrm R_{\pmin}(\Gamma)\colonequals\frac{\ip{\Gamma}{\Bfield}}{|\pmin^2\Gamma|}
$$
for every $\Gamma\in\mathcal{N}$, and set
$$
	\critfield\colonequals\frac{|\log\ep|}{2\mathrm R_\ep},\qquad \mathrm R_\ep\colonequals\sup_{\Gamma\in\mathcal{N}}\mathrm R_{\pmin}(\Gamma).
$$
We refer to \cite[Section 1.4]{diaz-roman-3d} for further details on normal currents and the weighted isoflux problem.

\medskip 
We now introduce the first quantitative regularity assumption. We assume that there exist $s\in(3,\infty)$, $M_s>0$, and $C>0$, independent of $\ep$, such that
\begin{equation}\label{assumptionLs}
	\|\nabla\pmin\|_{L^{s}(\Omega,\R^3)}\leq C|\log\ep|^{M_s}.
\end{equation}

Since $b\leq \pmin^2\leq 1$, it follows that
\begin{equation}\label{regularitypminW1s}
	\|\pmin\|_{W^{1,s}(\Omega)}\leq C|\log\ep|^{M_s}.
\end{equation}

By Morrey's inequality,
$$
W^{1,s}(\Omega)\hookrightarrow C^{0,1-3/s}(\overline\Omega),
$$
and therefore
$$
\|\pmin\|_{C^{0,1-3/s}(\overline\Omega)}
\leq C|\log\ep|^{M_s}.
$$
In particular, \eqref{regularitypminW1s} automatically yields the H\"older regularity required to apply the $\ep$-level estimates. The converse implication does not hold, so \eqref{regularitypminW1s} is a slightly stronger assumption on the pinning profile.

The motivation for introducing \eqref{regularitypminW1s} stems from a gap in the proof of \cite[Proposition 2.2]{diaz-roman-3d}, which provides the H\"older estimate for $\Bfield$ used in the proof of \cite[Theorem 1.2]{diaz-roman-3d}. As observed in \cite{rom2026metastab}, for a general lower bound $b\in(0,1)$, Gr\"oger's regularity theorem \cite{Groger} only yields a uniform bound for $\Bfield$ in $W^{1,\pzero}(\Omega,\R^3)$ for some exponent $\pzero>2$ depending on $b$. Since $\pzero$ need not exceed $3$, Sobolev embedding does not guarantee H\"older regularity of $\Bfield$ in general. In \cite{rom2026metastab}, uniform H\"older estimates were recovered under an additional small-BMO assumption on $\pmin^2$ (via \cite[Theorem 2]{SunWang}), which is automatically satisfied when $b$ is sufficiently close to $1$.

To treat general values of $b\in(0,1)$ without imposing a small-BMO condition, Section~\ref{sec:reg-B0} exploits the stronger regularity assumption \eqref{regularitypminW1s} to recover the H\"older regularity of $\Bfield$. More precisely, Proposition~\ref{prop:Bfield-holder} shows that there exists $\McurlBfield >0$, depending only on $M_s$ and $s$, such that for every $\gamma\in(0,1)$,
$$
\|\Bfield\|_{C^{0,\gamma}_T(\Omega,\R^3)}
\leq C_\gamma |\log\ep|^{\McurlBfield}.
$$
Hereafter, for $\gamma\in(0,1]$, we denote by $C_T^{0,\gamma}(\Omega,\R^3)$ the space of vector fields in $C^{0,\gamma}(\Omega,\R^3)$ satisfying $B\times\nu=0$ on $\partial\Omega$, endowed with the usual $C^{0,\gamma}$ norm.

Although this estimate is no longer uniform in $\ep$, the loss is only a power of $|\log\ep|$. Since the error term in the $\ep$-level vorticity estimate can be made to decay faster than any prescribed inverse power of $|\log\ep|$ by choosing the relevant parameters appropriately, this logarithmic loss is easily absorbed, and the proof of \cite[Theorem 1.2]{diaz-roman-3d} goes through essentially unchanged.

We are ready to state our main results. Our first result is a refinement of \cite[Theorem 1.2]{diaz-roman-3d} that establishes quantitative weak vortexlessness below the first critical field. In addition to resolving the gap in the estimate for $\Bfield$ under \eqref{regularitypminW1s}, it strengthens the original statement by providing decay faster than any prescribed inverse power of $|\log\ep|$. As we shall see, item~\ref{thm1:itemb} is essential for our clearing-out argument, since it allows us to bound the weighted free energy 
$$
\fenweightthreed(u,A)\colonequals\frac12\int_\Omega\left(\pmin^2|\nabla_Au|^2+\frac{\pmin^4}{2\ep^2}(1-|u|^2)^2\right)+\frac12\int_{\R^3}|\curl A|^2,
$$
by arbitrarily small inverse powers of $|\log\ep|$ once the parameters are chosen appropriately.

\begin{theorem}[Variant of Theorem~1.2 in \cite{diaz-roman-3d}]\label{thm:main}
	Suppose $\liminf_{\ep \to 0^+}\mathrm{R}_\ep>0$.
	
	Assume that \eqref{assumptionLs} holds and that there exist $C_2>0$ and $\kappa>0$, independent of $\ep$, such that $a_\ep$ is constant in the set
	$$
	\left\{x\in\Omega:\dist(x,\partial\Omega)\leq C_2|\log\ep|^{-\kappa}\right\}.
	$$
	
	Then for all $N\in\N$ and $\gamma_0\in(0,1)$ there exist $\ep_0,K_0>0$, depending on $N$ and $\gamma_0$, such that, for any $\ep<\ep_0$ and $h_{\mathrm{ex}}\le \critfield-K_0\log|\log\ep|$, the global minimizers $(\u,\A)$ of $\glthreed$ in $H^1(\Omega,\C)\times[A_{\mathrm{ex}}+H_{\curl}]$ are such that, letting
	$$
	({u},{A})=\big(\pmin^{-1}\u e^{-ih_{\mathrm{ex}}\phi_\ep^0},\A-h_{\mathrm{ex}}A^0_\ep\big),
	$$
	as $\ep\to0$, we have
	\begin{enumerate}[leftmargin=*,itemsep=3pt,label=\normalfont(\alph*)]
		\item\label{thm1:itema} $\|\mu({u},{A})\|_{(C^{0,\gamma}_T(\Omega,\R^3))^*}=o\left(|\log\ep|^{-N-1}\right)$ for any $\gamma\in[\gamma_0,1]$.
		\item\label{thm1:itemb} $\fenweightthreed({u},{A})=o\left(|\log\ep|^{-N}\right).$
		\item\label{thm1:itemc} $\glthreed(\u,\A)=\glthreed\big(\pmin e^{ih_{\mathrm{ex}}\phi^0_\ep},h_{\mathrm{ex}}A^0_\ep\big)+o\left(|\log\ep|^{-N}\right).$
	\end{enumerate}
\end{theorem}
\begin{remark}
	The regularity assumption in the statement differs from that of \cite[Theorem 1.2]{diaz-roman-3d} in that the quantitative H\"older estimate on $\pmin$ is replaced by \eqref{assumptionLs}. By Morrey's inequality, \eqref{assumptionLs} implies the H\"older regularity required in Theorem~\ref{epleveltools}.
\end{remark}

Theorem~\ref{thm:main} gives \emph{weak vortexlessness}: both the vorticity and the weighted free energy of the fluctuation around the Meissner state vanish in appropriate weak and energetic senses. Our main objective is to upgrade this information to a pointwise statement showing that the modulus of the order parameter remains uniformly close to the pinning profile.

For this purpose, we impose the stronger quantitative Lipschitz assumption
\begin{equation}\label{regularitypmin}
	\|\nabla\pmin\|_{L^\infty(\Omega,\R^3)}\leq C|\log\ep|^{M_\infty}
\end{equation}
for some $C,M_\infty>0$ independent of $\ep$. Recalling that $b\leq\pmin^2\leq1$, assumption \eqref{regularitypmin} implies \eqref{regularitypminW1s} for every $s\in (3,\infty]$. The stronger control in \eqref{regularitypmin} is used specifically in the weighted Pohozaev identity and the resulting pseudomonotonicity argument.

\begin{theorem}[Clearing-out up to the boundary for global minimizers]\label{thm:clearing-final}
	Assume the hypotheses of Theorem~\ref{thm:main} and \eqref{regularitypmin}. Let $M_\infty$ be as in \eqref{regularitypmin} and let $k\geq k_0\colonequals7M_\infty+1$ be an integer. Then there exist $\ep_0>0$ and $K_0=K_0(k)>0$ such that the following holds for $\ep<\ep_0$ and $h_\ex\leq\critfield-K_0\log|\log\ep|$. If $(\u,\A)$ is a global minimizer of $\glthreed$, letting
	$$
	(u,A)\colonequals\big(\pmin^{-1}\u e^{-ih_\ex\phi_\ep^0},\A-h_\ex A_\ep^0\big),
	$$
	then, as $\ep\to0$,
	\begin{equation}\label{eq:clearing-final}
		\|1-|u|\|_{L^\infty(\Omega)}=o(1)\qquad\text{and}\qquad\||\u|-\pmin\|_{L^\infty(\Omega)}=o(1).
	\end{equation}
	Moreover, for every $\gamma\in(0,1]$,
	\begin{equation}\label{eq:improved-conclusions}
		\fenweightthreed(u,A)=o(\ep),\qquad\|\mu(u,A)\|_{(C_T^{0,\gamma}(\Omega,\R^3))^*}=o(\ep).
	\end{equation}
	Here, $K_0$ is the constant furnished by Theorem~\ref{thm:main} with decay exponent $N=k+1$ and any fixed $\gamma_0\in(0,1)$.
	\end{theorem}

Theorem~\ref{thm:clearing-final} gives \emph{strong vortexlessness}: the modulus of the reduced order parameter $u$ converges uniformly to $1$, or equivalently, the modulus of the original order parameter $\u$ converges uniformly to the pinning profile $\pmin$. In particular, global minimizers have no zeros for sufficiently small $\ep$ when the applied field lies below $\critfield$ by a sufficiently large multiple of $\log|\log\ep|$.

Combining the lower bound furnished by Theorem~\ref{thm:main} with the upper bound obtained in \cite[Theorem 1.2]{RomanPinningUpper25}, we have that
\begin{equation*}
	H_{c_1}=\critfield+O(\log|\log\ep|).
\end{equation*}

Beyond the asymptotic characterization of $H_{c_1}$, Theorem~\ref{thm:clearing-final} yields a strong description of the Meissner regime below the first critical field by proving that the modulus of global minimizers remains uniformly close to the pinning profile. This completes, in the inhomogeneous three-dimensional setting, the counterpart of the picture established in \cite{Rom-CMP} for the homogeneous model.

\subsection*{Plan of the paper} The rest of the paper is organized as follows. In Section~\ref{sec:prelims}, we collect the preliminary material needed throughout the paper. In Section~\ref{sec:reg-B0}, we establish quantitative regularity estimates for the Meissner state approximation. Section~\ref{sec:toolsclearingout} is devoted to the tools required for the clearing-out argument. Finally, in Sections~\ref{sec:proofthm1} and \ref{sec:proofthm2}, we prove Theorems~\ref{thm:main} and \ref{thm:clearing-final}, respectively.

\subsection*{Acknowledgments} 

This work was partially funded by ANID FONDECYT 1231593.

\subsection*{AI declaration}

The mathematical content of this paper, including its results and proof strategies, was developed by the authors. Microsoft 365 Copilot, Google Gemini 4, and Claude Opus 5.5 were used to assist with the revision of the manuscript, including language editing, improvements in clarity and organization, and checks of notation, consistency, and selected mathematical calculations. The authors carefully reviewed and independently verified all AI-assisted changes and take full responsibility for the content of the manuscript.

\section{Preliminaries}\label{sec:prelims}
This section follows closely the lines of \cite{rom2026metastab}, in the sense that many of the definitions and proofs are taken directly from that work or slightly modified to suit the needs of the present paper.

We first present the three gauge choices that will be used throughout the paper. We then establish some properties and estimates for critical points of the Ginzburg--Landau energy. We conclude the section by recalling an essential energy splitting, the inhomogeneous $\ep$-level estimates that will allow us to prove Theorem~\ref{thm:main}, and a basic clearing-out estimate on balls of radius $\ep$ that is essential for the first step of our strategy.

\subsection{Gauge choices}

We shall make use of three gauge choices throughout the paper. The first is the classical Coulomb gauge.

\begin{definition}
	We say that a vector field $A$ is in the Coulomb gauge in $\Omega$ if
	\begin{equation*}
		\left\{
		\begin{array}{rcll}
			\diver A&=&0&\mathrm{in}\ \Omega,\\
			A\cdot\nu&=&0&\mathrm{on}\ \partial\Omega.
		\end{array}
		\right.
	\end{equation*}
\end{definition}
Every vector field $A\in L^2(\Omega,\R^3)$ can be transformed into this gauge by solving the corresponding Neumann problem. Extending the resulting phase to a function in $H^1_{\loc}(\R^3)$, the associated gauge transformation is defined in the whole space.

\medskip
The presence of the pinning profile $\pmin$ naturally leads to the following weighted version, which will play a central role in the proof of the clearing-out result.
\begin{definition}
	We say that $A$ is in the weighted Coulomb gauge associated with the weight $\pmin^2$ if
	\begin{equation}\label{weightedCoulombgauge}
		\left\{
		\begin{array}{rcll}
			\diver(\pmin^2A)&=&0&\mathrm{in}\ \Omega,\\
			A\cdot\nu&=&0&\mathrm{on}\ \partial\Omega.
		\end{array}
		\right.
	\end{equation}
\end{definition}
Since $b\leq\pmin^2\leq1$, every vector field $A\in L^2(\Omega,\R^3)$ can also be transformed into this gauge via the corresponding weighted Neumann problem. As before, the phase may be extended to an element of $H^1_{\loc}(\R^3)$, allowing the associated gauge transformation to be defined in the whole space.

\medskip
For the global gauge choice, we recall some standard properties of the homogeneous Sobolev space introduced in the previous section. By Sobolev embedding,
\begin{equation}\label{space1}
	\|A\|_{L^6(\R^3,\R^3)}\leq C\|\nabla A\|_{L^2(\R^3,\R^3)}
\end{equation}
for every $A\in\dot H^1(\R^3,\R^3)$. Moreover, by \cite[Proposition 2.4]{KozSoh},
$$
\dot H^1(\R^3,\R^3)=\left\{A\in L^6(\R^3,\R^3)\ \middle|\ \nabla A\in L^2(\R^3,\R^3)\right\}.
$$
Combining \eqref{space1} and \eqref{space2}, we obtain
$$
\|A\|_{L^6(\R^3,\R^3)}\leq C\|\curl A\|_{L^2(\R^3,\R^3)}
$$
for every $A\in\dot H^1_{\diver=0}$.

\begin{definition}
	We say that a vector field $A$ is in the global divergence-free gauge if
	$$
	A\in\dot H^1_{\diver=0}.
	$$
\end{definition}
The following decomposition, proved in \cite[Lemma 2.1]{Rom-CMP}, shows that every admissible configuration may be transformed into this gauge.
\begin{lemmaA}\label{lemma:decompositionR3}
	Every vector field $A\in H_{\curl}$ admits a decomposition
	$$
	A=\curl B+\nabla\Phi,
	$$
	where $B,\curl B\in\dot H^1_{\diver=0}$ and $\Phi\in H^2_{\loc}(\R^3)$.
\end{lemmaA}
Applying the gauge transformation
$$
u\to ue^{-i\Phi},\qquad A\to A-\nabla\Phi,
$$
we obtain a gauge-equivalent vector potential belonging to $\dot H^1_{\diver=0}$.

\subsection{Critical points of the Ginzburg--Landau energy}

We now present some properties of critical points of $\glthreed$ that will be useful throughout the paper. The Euler--Lagrange equations associated with $\glthreed$ are given in the following proposition.

\begin{propA}[Ginzburg--Landau equations; {\cite[Proposition 2.1]{rom2026metastab}}]
	If the configuration $(\u,\A)\in H^1(\Omega,\C)\times[A_{\ex}+H_{\curl}]$ is a critical point of $\glthreed$, then $(\u,\A)$ satisfies the system
	\begin{equation}\label{GLeq}
		\left\lbrace
		\begin{array}{rcll}
			-(\nabla_\A)^2\u&=&\displaystyle\frac1{\ep^2}\u(a_\ep-|\u|^2)&\mathrm{in}\ \Omega\\
			\curl(\H-\H_\ex)&=&\ip{i\u}{\nabla_\A\u}\chi_\Omega&\mathrm{in}\ \R^3\\
			\nabla_\A\u\cdot\nu&=&0&\mathrm{on}\ \partial\Omega\\
			\lbrack\H-\H_\ex\rbrack\times\nu&=&0&\mathrm{on}\ \partial\Omega.
		\end{array}
		\right.
	\end{equation}
	Here, $\chi_\Omega$ is the characteristic function of $\Omega$, $[\,\cdot\,]$ denotes the jump across $\partial\Omega$,
	$$
	\nabla_\A\u\cdot\nu=\sum_{j=1}^3(\partial_j\u-i\A_j\u)\nu_j,
	$$
	and the covariant Laplacian $(\nabla_\A)^2$ is defined by
	$$
	(\nabla_\A)^2\u=(\diver-i\A\cdot)\nabla_\A\u.
	$$
\end{propA}

\begin{remark}
	Let $(\u,\A)$ be a critical point of $\glthreed$ and set $u=\pmin^{-1}\u$. Since $\pmin$ solves \eqref{pde rho}, the first equation and the first boundary condition in \eqref{GLeq} can be rewritten as
	\begin{equation}\label{eq:post_splitting_eq}
		\left\lbrace
		\begin{array}{rcll}
			-(\diver-i\A\cdot)(\pmin^2\nabla_\A u)&=&\displaystyle\frac1{\ep^2}\pmin^4u(1-|u|^2)&\mathrm{in}\ \Omega\\
			\nabla_\A u\cdot\nu&=&0&\mathrm{on}\ \partial\Omega.
		\end{array}
		\right.
	\end{equation}
	This system will serve as the starting point for the Pohozaev identity.
\end{remark}

We next recall bounds and regularity estimates for the order parameter. They will be used in the derivation of the monotonicity formula and the clearing-out estimate.

\begin{propA}[{\cite[Proposition 2.3]{rom2026metastab}}]\label{prop:regularity}
	Let $(\u,\A)$ be a critical point of $\glthreed$. Then
	\begin{equation}\label{boundu}
		\|\u\|_{L^\infty(\Omega,\C)}\leq1,\quad\mbox{}\quad\left\|\frac{\u}{\pmin}\right\|_{L^\infty(\Omega,\C)}\leq1.
	\end{equation}
	Moreover, if $\A$ is in the Coulomb gauge in $\Omega$, then
	\begin{equation*}
		\|\nabla\u\|_{L^\infty(\Omega,\C^3)}\leq C\left(\frac1\ep+\|\A\|_{L^\infty(\Omega,\R^3)}\right).
	\end{equation*}
\end{propA}

The following result provides a regularity estimate for the magnetic potential $\A$, which will be essential in the derivation of the pseudomonotonicity formula. It is formulated slightly differently from \cite[Proposition 2.4]{rom2026metastab}, but it is obtained directly in its proof.

\begin{propA}[{\cite[Proposition 2.4]{rom2026metastab}}]\label{prop:regularityA}
	Let $(\u,\A)$ be a critical point of $\glthreed$. If $\A$ is in the Coulomb gauge in $\Omega$, then
	$$
	\|\A\|_{W^{1,6}(\Omega,\R^3)}\leq C\left(\left\|\nabla_\A\left(\frac{\u}{\pmin}\right)\right\|_{L^2(\Omega,\C^3)}+h_\ex\right).
	$$
\end{propA}

\subsection{Energy estimates}

We conclude this section by presenting three ingredients that will play a central role in the proofs of our main results. We first recall the energy splitting established in \cite[Proposition 1.1]{diaz-roman-3d}. We then present the inhomogeneous $\ep$-level estimates from \cite[Theorem 1.1]{diaz-roman-3d}, which will be used to bound the free energy. Finally, we recall a basic clearing-out estimate on small balls from \cite[Lemma 2.1]{rom2026metastab}.

\begin{propA}[Energy splitting; {\cite[Proposition 1.1]{diaz-roman-3d}}]
	Given any configuration
	$$
	(\u,\A)\in\minsp,
	$$
	let $(u,A)$ be defined by
	$$
	(\u,\A)=(\pmin u\phase,A+h_\ex A_\ep^0).
	$$
	Then
	\begin{equation}\label{Split3D}
		\glthreed(\u,\A)=\glthreed(\pmin\phase,h_\ex A_\ep^0)+\fenweightthreed(u,A)-h_\ex\int_\Omega\mu(u,A)\cdot\Bfield+R_0,
	\end{equation}
	where
	$$
	\fenweightthreed(u,A)=\frac12\int_\Omega\left(\pmin^2|\nabla_Au|^2+\frac{\pmin^4}{2\ep^2}(1-|u|^2)^2\right)+\frac12\int_{\R^3}|\curl A|^2
	$$
	and
	$$
	R_0\colonequals\frac{h_\ex^2}{2}\int_\Omega\frac{|\curl\Bfield|^2}{\pmin^2}(|u|^2-1).
	$$
\end{propA}

\begin{theoremA}[Inhomogeneous $\ep$-level estimates; {\cite[Theorem 1.1]{diaz-roman-3d}}]\label{epleveltools}
	Assume that $a_\ep$ satisfies the following conditions:
	\begin{enumerate}[leftmargin=*,itemsep=3pt,label=\normalfont(\roman*)]
		\item There exist $\alpha\in(0,1)$, $N>0$, and $C_1>0$, independent of $\ep$, such that
		$$
		\|\pmin\|_{C^{0,\alpha}(\Omega)}\leq C_1|\log\ep|^N.
		$$
		\item There exist $C_2>0$ and $\kappa>0$, independent of $\ep$, such that $a_\ep$ is constant in the set
		$$
		\{x\in\Omega\mid\dist(x,\partial\Omega)\leq C_2|\log\ep|^{-\kappa}\}.
		$$
	\end{enumerate}
	Then the following holds. For any $m,n,\tilde{C}>0$, there exist $C,\ep_0>0$, depending only on $\tilde{C},m,n$, and $\partial\Omega$, such that, for any $\ep<\ep_0$, if $(u,A)\in H^1(\Omega,\C)\times H^1(\Omega,\R^3)$ satisfies
	$$
	\fenweightthreed(u,A)\leq\tilde{C}|\log\ep|^m,
	$$
	then there exists a polyhedral $1$-current $\nu_\ep$ such that:
	\begin{enumerate}[leftmargin=*,itemsep=3pt,label=\normalfont(\arabic*)]
		\item $\nu_\ep/\pi$ has integer multiplicity,
		\item $\partial\nu_\ep=0$ relative to $\Omega$,
		\item $\mathrm{supp}(\nu_\ep)\subset S_{\nu_\ep}\subset\overline\Omega$, with $|S_{\nu_\ep}|\leq C|\log\ep|^{-q}$, where $q=q(m,n)\colonequals3(m+n)$,
		\item
		\begin{equation}\label{weighted lower bound}
			\int_{S_{\nu_\ep}}\left(\pmin^2|\nabla_Au|^2+\frac{\pmin^4}{2\ep^2}(1-|u|^2)^2+|\curl A|^2\right)\geq|\pmin^2\nu_\ep|\left(|\log\ep|-C\log|\log\ep|\right)-\frac{C}{|\log\ep|^n},
		\end{equation}
		where $|\pmin^2\nu_\ep|$ denotes the mass of the $1$-current $\pmin^2\nu_\ep$ in $\Omega$,
		\item for any $\gamma\in(0,1]$, there exists $C_\gamma>0$, depending only on $\gamma$ and $\partial\Omega$, such that
		\begin{equation}\label{vorticity estimate}
			\norm{\mu(u,A)-\nu_\ep}{(C^{0,\gamma}_T(\Omega,\R^3))^*}\leq C_\gamma\frac{\fenweightthreed(u,A)+1}{|\log\ep|^{q\gamma}}.
		\end{equation}
	\end{enumerate}
\end{theoremA}

\begin{lemmaA}[Basic clearing-out estimate on small balls; {\cite[Lemma 2.1]{rom2026metastab}}]\label{lemma:clearingout}
	Let $u\in H^1(\Omega,\C)$ satisfy
	$$
	\|\nabla|u|\|_{L^\infty(\Omega,\R^3)}\leq C\ep^{-1}.
	$$
	Then
	$$
	(1-|u|(x_0))^5\leq\frac{C}{\ep^3}\int_{B_\ep(x_0)\cap\Omega}(1-|u|^2)^2,\quad\forall x_0\in\overline\Omega.
	$$
\end{lemmaA}

\section{Regularity estimates for the Meissner state approximation}\label{sec:reg-B0}
This section is devoted to regularity estimates for the Meissner state approximation $(\meisconf)$ and the associated vector field $\Bfield$ introduced in \eqref{defBfield}--\eqref{conditionsB0}. We first recall the uniform estimates established in \cite{rom2026metastab}. For a general lower bound $b\in(0,1)$ for $a_\ep$, these estimates only provide a bound for $\Bfield$ in $W^{1,\pzero}(\Omega,\R^3)$ for some $\pzero>2$, which need not yield H\"older control. Under the quantitative Sobolev assumption \eqref{regularitypminW1s}, we then derive estimates with an explicit exponent larger than $3$, at the expense of a power of $|\log\ep|$. These estimates recover the H\"older regularity of $\Bfield$ required in the proof of Theorem~\ref{thm:main} and provide the logarithmic-loss exponents needed for Theorem~\ref{thm:clearing-final}.

\subsection{Uniform estimates in \texorpdfstring{$\ep$}{ε} for the Meissner state approximation}

We begin by recalling the estimates for the Meissner state approximation obtained in \cite{rom2026metastab}. Recall that $A_\ep^0$ denotes the minimizer of the reduced functional \eqref{reducedenergy} satisfying \eqref{spaceforA0}, and that $\phi_\ep^0\colonequals\phi_{A_\ep^0}$. Assuming that $A_{0,\ex}$ is in the Coulomb gauge and that $\curl A_{0,\ex}=H_{0,\ex}\in L^\infty(\R^3,\R^3)$, we have $A_{0,\ex}\in W^{1,6}(\Omega,\R^3)$, and hence
\begin{equation}\label{estimateforA0}
	\|A_\ep^0\|_{W^{1,6}(\Omega,\R^3)}+\|A_\ep^0\|_{L^\infty(\Omega,\R^3)}+\|\nabla\phi_\ep^0\|_{L^2(\Omega,\R^3)}+\|\curl\Bfield\|_{L^2(\Omega,\R^3)}\leq C,
\end{equation}
where $C=C(\Omega,b)>0$ is independent of $\ep$. Indeed, the minimizing property of $A_\ep^0$, together with \eqref{space2}, gives
$$
\|A_\ep^0-A_{0,\ex}\|_{\dot H^1_{\diver=0}}+\|\nabla\phi_\ep^0\|_{L^2(\Omega,\R^3)}+\|\curl\Bfield\|_{L^2(\Omega,\R^3)}\leq C,
$$
while elliptic regularity and Sobolev embedding yield the first two terms in \eqref{estimateforA0}.

Moreover, since $b\leq\pmin^2\leq1$, Gr\"oger's result \cite{Groger} provides an exponent $\pzero>2$, depending only on $b$, such that
\begin{equation}\label{estimateforgradphip0}
	\|\nabla\phi_\ep^0\|_{L^{\pzero}(\Omega,\R^3)}\leq C\|A_\ep^0\|_{L^{\pzero}(\Omega,\R^3)},
\end{equation}
where both $\pzero$ and $C$ are independent of $\ep$. Combining \eqref{defBfield}, \eqref{estimateforA0}, and \eqref{estimateforgradphip0}, and using \cite[Corollary 3.4]{amrouche-seloula} together with \eqref{conditionsB0}, we obtain
$$
\|\curl\Bfield\|_{L^{\pzero}(\Omega,\R^3)}+\|\Bfield\|_{W^{1,\pzero}(\Omega,\R^3)}\leq C.
$$
In particular, if $\pzero>3$, then Sobolev embedding gives
\begin{equation}\label{estimateHolderB0}
	\|\Bfield\|_{C^{0,\alpha}_T(\Omega,\R^3)}\leq C
\end{equation}
for every $\alpha\leq1-3/\pzero$.

More generally, the result of Byun and Wang \cite[Theorem 2]{SunWang} allows one to obtain \eqref{estimateforgradphip0} for any prescribed finite exponent $\pzero>2$, provided that $\pmin^2$ satisfies a suitable small-BMO condition whose smallness threshold depends on $\pzero$. As observed in \cite{rom2026metastab}, for each prescribed finite $\pzero>2$, this condition is satisfied when $b$ is sufficiently close to $1$, with the required closeness depending on $\pzero$. In particular, one may choose $\pzero>3$ and obtain the uniform H\"older estimate \eqref{estimateHolderB0}.

For a general $b\in(0,1)$, however, Gr\"oger's result only provides some exponent $\pzero>2$, which may be close to $2$, while the small-BMO condition required to prescribe an exponent larger than $3$ is not guaranteed. We therefore derive below estimates with explicit exponents larger than $3$ under \eqref{regularitypminW1s}. These estimates are no longer uniform in $\ep$, but their dependence on $\ep$ is controlled by powers of $|\log\ep|$.

\subsection{Estimates for gauge transformations}

We now derive an estimate for the gauge transformation into the weighted Coulomb gauge associated with the weight $\pmin^2$. This estimate will later be applied to $\phi_\ep^0$ in order to obtain explicit regularity estimates for $\Bfield$.

We first recall a consequence of \cite[Theorem~3.17(ii)]{Tro}, for $p\geq2$, and \cite[Problem~3.5]{Tro}, for $1<p<2$.\footnote{The Sobolev spaces denoted by $W^{k,p}$ in this paper are denoted by $H^{k,p}$ in \cite{Tro}.} In the notation adapted to this paper, the natural boundary condition is imposed on $\Gamma=\partial\Omega$. Consequently,
$$
W_0^{1,q'}(\Omega\cup\Gamma)=W^{1,q'}(\Omega),
$$
since $W_0^{1,q'}(\Omega\cup\Gamma)$ consists of the functions in $W^{1,q'}(\Omega)$ whose trace vanishes on $\partial\Omega\setminus\Gamma=\varnothing$.

\begin{theoremA}[$W^{2,p}$ regularity for the Neumann problem; {\cite[Theorem~3.17(ii) and Problem~3.5]{Tro}}]\label{thm:H2p-neumann}
	Let $\Omega\subset\R^d$ be a bounded domain with boundary of class $C^{1,1}$, let $1<p<\infty$, set $q\colonequals\min\{p,2\}$, and let $q'$ be its conjugate exponent. Suppose that $u\in W^{1,q}(\Omega)$ satisfies
	$$
	\int_\Omega\nabla u\cdot\nabla v=\int_\Omega f^0v+\int_\Omega F\cdot\nabla v\qquad\text{for every }v\in W^{1,q'}(\Omega),
	$$
	where $f^0\in L^p(\Omega)$ and $F\in W^{1,p}(\Omega,\R^d)$. Then $u\in W^{2,p}(\Omega)$ and
	$$
	\|u\|_{W^{2,p}(\Omega)}\leq C\left(\|f^0\|_{L^p(\Omega)}+\|F\|_{W^{1,p}(\Omega,\R^d)}+\|u\|_{W^{1,q}(\Omega)}\right),
	$$
	where $C=C(\Omega,p)>0$ is independent of $u$, $f^0$, and $F$.
\end{theoremA}

For $1<p<2$, this result appears as \cite[Problem~3.5]{Tro}, with the convention introduced in \cite[Problem~3.3]{Tro} that the $H^1$ norm of $u$ in the estimate is replaced by its $W^{1,p}$ norm. The proof follows the argument for $p\geq2$ by duality: the tangential second derivatives are identified with solutions of Neumann-type problems with data in $L^p$, and the normal second derivative is then recovered from the equation.

We apply this result to the equation defining the gauge transformation into the weighted Coulomb gauge.

\begin{lemma}[$W^{2,s}$ regularity for the weighted Neumann problem]\label{lem:fase-Ws}
	Let $3<s<\infty$, $\pmin\in W^{1,s}(\Omega)$, and $A\in W^{1,s}(\Omega,\R^3)$. There exists a unique $\varphi\in H^1(\Omega)$ with $\int_\Omega\varphi=0$ such that
	\begin{equation}\label{eq:fase-debil}
		\int_\Omega\pmin^2\nabla\varphi\cdot\nabla v=\int_\Omega\pmin^2A\cdot\nabla v\qquad\text{for all }v\in H^1(\Omega).
	\end{equation}
	Moreover, $\varphi\in W^{2,s}(\Omega)$, $\nabla\varphi\in C^{0,1-3/s}(\overline\Omega,\R^3)$, and
	\begin{equation}\label{eq:fase-cota}
		\|\varphi\|_{W^{2,s}(\Omega)}\leq C\big(1+\|\nabla\pmin\|_{L^s(\Omega)}\big)^m\|A\|_{W^{1,s}(\Omega,\R^3)},\qquad m=m(s)\colonequals\left\lfloor\frac{3s}{s-3}\right\rfloor+1,
	\end{equation}
	where $C=C(\Omega,b,s)>0$ is independent of $\pmin$ and $A$.
\end{lemma}

\begin{proof}
	Throughout the proof, $C>0$ denotes a constant depending only on $\Omega$, $b$, and $s$, whose value may change from line to line. Since $s>3$, Morrey's inequality \cite[Theorem~7.26(ii)]{GilbargTrudinger01} gives
	$$
	\|A\|_{L^\infty(\Omega,\R^3)}\leq C\|A\|_{W^{1,s}(\Omega,\R^3)}.
	$$
	Set
	$$
	K\colonequals\left\|\frac{\nabla\pmin^2}{\pmin^2}\right\|_{L^s(\Omega)}=\left\|\frac{2\nabla\pmin}{\pmin}\right\|_{L^s(\Omega)}\leq2b^{-1/2}\|\nabla\pmin\|_{L^s(\Omega)}.
	$$
	\begin{enumerate}[label=\textnormal{\bfseries Step~\arabic*.},ref=\arabic*,leftmargin=0pt,labelsep=*,itemindent=*,itemsep=10pt,topsep=10pt]
		\item \textbf{Existence and energy estimate.} Existence and uniqueness follow from the Lax--Milgram theorem on
		$$
		V\colonequals\left\{v\in H^1(\Omega):\int_\Omega v=0\right\},
		$$
		using $b\leq\pmin^2\leq1$ and the Poincar\'e--Wirtinger inequality. Since both sides of \eqref{eq:fase-debil} vanish on constant functions, the resulting identity holds for every $v\in H^1(\Omega)$. Testing \eqref{eq:fase-debil} with $v=\varphi$ gives
		$$
		\|\pmin\nabla\varphi\|_{L^2}^2\leq\|\pmin A\|_{L^2}\|\pmin\nabla\varphi\|_{L^2},
		$$
		and hence, using the Poincar\'e--Wirtinger inequality,
		\begin{equation}\label{eq:fase-energia}
			\|\varphi\|_{H^1(\Omega)}\leq C\|A\|_{L^2(\Omega,\R^3)}\leq C\|A\|_{W^{1,s}(\Omega,\R^3)}.
		\end{equation}
		\item \textbf{Reduction to a Neumann problem for the Laplacian.} For every $w\in H^1(\Omega)$, we have $w\pmin^{-2}\in H^1(\Omega)$ and
		$$
		\nabla(w\pmin^{-2})=\pmin^{-2}\nabla w+w\nabla\pmin^{-2}.
		$$
		Indeed, $w\in L^6(\Omega)$ and $\nabla\pmin^{-2}\in L^s(\Omega)\subset L^3(\Omega)$, so
		$$
		\|w\nabla\pmin^{-2}\|_{L^2}\leq\|w\|_{L^6}\|\nabla\pmin^{-2}\|_{L^3}<\infty.
		$$
		Testing \eqref{eq:fase-debil} with $v=w\pmin^{-2}$ gives
		$$
		\int_\Omega\nabla\varphi\cdot\nabla w=\int_\Omega hw+\int_\Omega A\cdot\nabla w\qquad\text{for all }w\in H^1(\Omega),
		$$
		where
		$$
		h\colonequals\frac{\nabla\pmin^2}{\pmin^2}\cdot(\nabla\varphi-A).
		$$
		If $\nabla\varphi\in L^q(\Omega)$ for some $q\in[2,\infty]$ and $1/r=1/s+1/q$, then
		\begin{equation}\label{eq:fase-h}
			\|h\|_{L^r(\Omega)}\leq K\left(\|\nabla\varphi\|_{L^q(\Omega)}+C\|A\|_{W^{1,s}(\Omega,\R^3)}\right).
		\end{equation}
		Consequently, if $h\in L^p(\Omega)$ and $A\in W^{1,p}(\Omega,\R^3)$ for some $p\in(1,\infty)$, Theorem~\ref{thm:H2p-neumann} gives
		\begin{equation}\label{eq:fase-R}
			\|\varphi\|_{W^{2,p}(\Omega)}\leq C(\Omega,p)\left(\|h\|_{L^p(\Omega)}+\|A\|_{W^{1,p}(\Omega,\R^3)}+\|\varphi\|_{W^{1,\min\{p,2\}}(\Omega)}\right).
		\end{equation}
		The last term is controlled by \eqref{eq:fase-energia}.
		\item \textbf{Bootstrap up to $\nabla\varphi\in L^\infty(\Omega)$.} Set
		\begin{gather*}
			\delta\colonequals\frac13-\frac1s\in\left(0,\frac13\right),\qquad t_k\colonequals\frac12-\frac{k\delta}{2},\qquad \frac1{r_k}\colonequals\frac1s+t_k,\\
			N\colonequals\min\{k\in\N:t_k<\delta\}=\left\lfloor\frac1\delta\right\rfloor-1.
		\end{gather*}
		For $0\leq k\leq N$, we have $t_k>0$, $r_k<s$, and
		$$
		\frac1{r_k}-\frac13=t_k-\delta.
		$$
		In particular, $r_k\leq3$ for $k<N$, whereas $r_N>3$. Define
		$$
		Q_k\colonequals\|\nabla\varphi\|_{L^{1/t_k}(\Omega)}.
		$$
		We claim that
		\begin{equation}\label{eq:fase-induccion}
			Q_k\leq C(1+K)^k\|A\|_{W^{1,s}(\Omega,\R^3)},\qquad 0\leq k\leq N.
		\end{equation}
		For $k=0$, this follows from \eqref{eq:fase-energia}. Suppose that \eqref{eq:fase-induccion} holds for some $k\leq N$. By \eqref{eq:fase-h}, $h\in L^{r_k}(\Omega)$ and
		$$
		\|h\|_{L^{r_k}(\Omega)}\leq K\left(Q_k+C\|A\|_{W^{1,s}(\Omega,\R^3)}\right).
		$$
		Since $r_k<s$, \eqref{eq:fase-R} and \eqref{eq:fase-energia} imply
		\begin{align*}
			\|\varphi\|_{W^{2,r_k}(\Omega)}&\leq C(1+K)\left(Q_k+\|A\|_{W^{1,s}(\Omega,\R^3)}\right)\\
			&\leq C(1+K)^{k+1}\|A\|_{W^{1,s}(\Omega,\R^3)}.
		\end{align*}
		If $k<N$, then $r_k\leq3$ and
		$$
		\frac1{r_k}-\frac13=t_k-\delta\leq t_{k+1}.
		$$
		Therefore, the Sobolev embedding gives
		$$
		Q_{k+1}=\|\nabla\varphi\|_{L^{1/t_{k+1}}(\Omega)}\leq C\|\varphi\|_{W^{2,r_k}(\Omega)}\leq C(1+K)^{k+1}\|A\|_{W^{1,s}(\Omega,\R^3)},
		$$
		which proves \eqref{eq:fase-induccion}. For $k=N$, since $r_N>3$, Morrey's inequality and the preceding estimate yield
		\begin{equation}\label{eq:fase-Linfty}
			\|\nabla\varphi\|_{L^\infty(\Omega)}\leq C\|\varphi\|_{W^{2,r_N}(\Omega)}\leq C(1+K)^{N+1}\|A\|_{W^{1,s}(\Omega,\R^3)}.
		\end{equation}
		\item \textbf{The $W^{2,s}$ estimate.} By \eqref{eq:fase-h} with $q=\infty$ and \eqref{eq:fase-Linfty},
		$$
		\|h\|_{L^s(\Omega)}\leq C(1+K)^{N+2}\|A\|_{W^{1,s}(\Omega,\R^3)}.
		$$
		Applying \eqref{eq:fase-R} with $p=s$ and using \eqref{eq:fase-energia}, we obtain
		$$
		\|\varphi\|_{W^{2,s}(\Omega)}\leq C(1+K)^{N+2}\|A\|_{W^{1,s}(\Omega,\R^3)}.
		$$
		Since $1/\delta=3s/(s-3)$, we have
		$$
		N+2=\left\lfloor\frac{3s}{s-3}\right\rfloor+1=m.
		$$
		The estimate \eqref{eq:fase-cota} now follows from
		$$
		1+K\leq C(b)\big(1+\|\nabla\pmin\|_{L^s(\Omega)}\big).
		$$
		Finally, Morrey's inequality gives $\nabla\varphi\in C^{0,1-3/s}(\overline\Omega,\R^3)$.
	\end{enumerate}
\end{proof}

\subsection{Estimates for \texorpdfstring{$\Bfield$}{B} with explicit exponents}

We now apply Lemma~\ref{lem:fase-Ws} to the gauge transformation defining the Meissner state approximation $(\meisconf)$ through \eqref{eqphi0}--\eqref{defBfield}, under assumption \eqref{regularitypminW1s}.

\begin{prop}\label{prop:Bfield-holder}
	Assume \eqref{regularitypminW1s} for some $s\in(3,\infty)$ and $M_s>0$. Then there exists $\McurlBfield>0$, independent of $\ep$, such that, for every $\gamma\in(0,1)$, there exists $C_\gamma>0$, independent of $\ep$, satisfying
	\begin{equation}\label{estimateforB0holder}
		\|\Bfield\|_{C^{0,\gamma}_T(\Omega,\R^3)}\leq C_\gamma|\log\ep|^{\McurlBfield}.
	\end{equation}
	The exponent $\McurlBfield$ depends only on $M_s$ and $s$.
\end{prop}

\begin{proof}
	Set $\widetilde s\colonequals\min\{s,6\}$. Since $\widetilde s\leq6$, \eqref{estimateforA0} gives
	$$
	\|A_\ep^0\|_{W^{1,\widetilde s}(\Omega,\R^3)}\leq C.
	$$
	By \eqref{regularitypminW1s} and Lemma~\ref{lem:fase-Ws}, with $m=m(\widetilde s)$, we obtain
	\begin{equation}\label{phi0estimatelep}
		\|\phi_\ep^0\|_{W^{2,\widetilde s}(\Omega)}\leq C\big(1+\|\nabla\pmin\|_{L^{\widetilde s}(\Omega)}\big)^m\|A_\ep^0\|_{W^{1,\widetilde s}(\Omega,\R^3)}\leq C|\log\ep|^{mM_s}.
	\end{equation}
	The product rule applied to \eqref{defBfield} yields
	$$
	|\nabla\curl\Bfield|\leq2|\nabla\pmin|\,|A_\ep^0-\nabla\phi_\ep^0|+\pmin^2|\nabla(A_\ep^0-\nabla\phi_\ep^0)|.
	$$
	Using \eqref{regularitypminW1s}, \eqref{phi0estimatelep}, and Morrey's inequality, we obtain
	\begin{align*}
		\|\curl\Bfield\|_{W^{1,\widetilde s}(\Omega,\R^3)}&\leq C\left(\|\nabla\pmin\|_{L^{\widetilde s}(\Omega,\R^3)}\|A_\ep^0-\nabla\phi_\ep^0\|_{L^\infty(\Omega,\R^3)}+\|A_\ep^0-\nabla\phi_\ep^0\|_{W^{1,\widetilde s}(\Omega,\R^3)}\right)\\
		&\leq C\left(|\log\ep|^{M_s+mM_s}+|\log\ep|^{mM_s}\right)\\
		&\leq C|\log\ep|^{(m+1)M_s}.
	\end{align*}
	Since $\widetilde s>3$, Morrey's inequality gives
	\begin{equation}\label{estimateforcurlB0Linf}
		\|\curl\Bfield\|_{L^\infty(\Omega,\R^3)}\leq C|\log\ep|^{(m+1)M_s}.
	\end{equation}
	Using the Friedrichs--Poincar\'e estimate as in \cite[Corollary 3.4]{amrouche-seloula}, we deduce, for every $p\in(1,\infty)$,
	$$
	\|\Bfield\|_{W^{1,p}(\Omega,\R^3)}\leq C_p|\log\ep|^{(m+1)M_s}.
	$$
	Set $\McurlBfield\colonequals(m+1)M_s$. Given $\gamma\in(0,1)$, choose $p\geq3/(1-\gamma)$. Then $\gamma\leq1-3/p$, and Morrey's inequality yields \eqref{estimateforB0holder}.
\end{proof}

\begin{remark}
	The preceding proposition yields H\"older estimates for $\Bfield$ with every exponent $\gamma\in(0,1)$, but does not provide a Lipschitz estimate, since the Friedrichs--Poincar\'e inequality used in its proof only yields $W^{1,p}$ estimates for finite $p$. One could instead try to derive $W^{2,p}$ estimates for $\Bfield$ from the $W^{1,p}$ estimates for $\curl\Bfield$. For instance, \cite[Corollaries 3.5 and 5.3]{amrouche-seloula} provide suitable higher-order div--curl estimates under the stronger assumption that $\partial\Omega$ is of class $C^{2,1}$. Choosing $p>3$ would then yield a Lipschitz estimate for $\Bfield$ by Sobolev embedding. Since the H\"older estimates above are sufficient for our purposes, we do not impose this additional boundary regularity.
\end{remark}

\section{Clearing-out tools}\label{sec:toolsclearingout}
In this section, we establish the tools needed to transport energy information from larger scales down to balls of radius $\ep$, following closely the strategy of \cite{chiron2005boundary}. The main result is a pseudomonotonicity formula for the normalized energy $\Phi$ defined in \eqref{eq:definitions_monotonicity}. The derivative of $\Phi$ involves energy terms in the interior of the ball and on its boundary. To control them, we first derive a Pohozaev identity, then estimate the error terms arising in the interior, and finally use the same identity to control the terms on $\partial\Omega$. Combining these estimates, we show that a suitable modification of $\Phi$, obtained by adding an error term and multiplying by an integrating factor, is nondecreasing. This allows us to control the normalized energy at scale $\ep$ in terms of the energy at larger scales.

Throughout this section, we assume the quantitative Lipschitz bound \eqref{regularitypmin} for the pinning profile $\pmin$. This assumption is used to control the error terms arising in the Pohozaev identity. Other parts of the paper require only weaker regularity assumptions on $\pmin$. We begin with some geometric properties of the $C^2$ boundary $\partial\Omega$.

\begin{lemma}[Geometry of the boundary]\label{lem:geom-borde}
	Let $\Omega\subset\R^3$ be a bounded domain with $\partial\Omega\in C^2$ and outward unit normal $\nu$. There exists $R>0$, depending only on $\Omega$, such that, setting
	$$
	\kappa\colonequals\frac1{2R},\qquad r_1\colonequals\min\{1,R/18,\operatorname{diam}\Omega\}.
	$$
	the following holds. Let $x_0\in\overline\Omega$, $d_0\colonequals\dist(x_0,\partial\Omega)$, let $p\in\partial\Omega$ satisfy $|x_0-p|=d_0$, and let $d_0\leq r<r_1$. If
	$$
	z\colonequals p-\frac r2\nu(p),
	$$
	then
	\begin{equation}\label{eq:geom-a}
		|x-z|\leq3r\qquad\text{for every }x\in B_r(x_0),
	\end{equation}
	and
	\begin{equation}\label{eq:geom-b}
		|(x-x_0)\cdot\nu(x)|\leq d_0+4\kappa r^2,\qquad (x-z)\cdot\nu(x)\geq\frac r4\qquad\text{for every }x\in\partial\Omega\cap B_r(x_0).
	\end{equation}
\end{lemma}

\begin{proof}
	We use the uniform two-sided sphere condition for bounded $C^2$ domains, as in \cite[Problem 2.11]{GilbargTrudinger01}: there exists $R>0$ such that, for every $q\in\partial\Omega$,
	\begin{equation}\label{eq:spheres}
		B_R(q-R\nu(q))\subset\Omega,\qquad B_R(q+R\nu(q))\subset\R^3\setminus\overline\Omega.
	\end{equation}
	\begin{enumerate}[label=\textnormal{\bfseries Step~\arabic*.},ref=\arabic*,leftmargin=0pt,labelsep=*,itemindent=*,itemsep=10pt,topsep=10pt]
		\item\label{step1lemgeom} \textbf{A quadratic estimate for the normal component.} We claim that
		$$
		|(x-p)\cdot\nu(x)|\leq\frac{|x-p|^2}{2R}\qquad\text{for every }x,p\in\partial\Omega.
		$$
		Indeed, from \eqref{eq:spheres} $p$ belongs neither to $B_R(x-R\nu(x))$ nor to $B_R(x+R\nu(x))$. Hence
		$$
		|p-x\pm R\nu(x)|\geq R,
		$$
		and therefore
		$$
		|p-x|^2\pm2R(p-x)\cdot\nu(x)\geq0,
		$$
		which proves the claim.
		\item \textbf{Proof of \eqref{eq:geom-a} and of the first estimate in \eqref{eq:geom-b}.} If $x\in B_r(x_0)$, then
		$$
		|x-z|\leq|x-x_0|+|x_0-p|+|p-z|<r+d_0+\frac r2\leq\frac52r\leq3r,
		$$
		which proves \eqref{eq:geom-a}. If $x\in\partial\Omega\cap B_r(x_0)$, then $|x-p|<r+d_0\leq2r$. Step~\ref{step1lemgeom} gives
		$$
		|(x-p)\cdot\nu(x)|\leq\frac{4r^2}{2R}=4\kappa r^2,
		$$
		and hence
		$$
		|(x-x_0)\cdot\nu(x)|\leq|(x-p)\cdot\nu(x)|+|p-x_0|\leq4\kappa r^2+d_0.
		$$
		\item \textbf{Positivity of the boundary contribution.} We prove the second estimate in \eqref{eq:geom-b}. Since $r/2\leq R$, the ball $B_{r/2}(z)$ is contained in $B_R(p-R\nu(p))\subset\Omega$. On the other hand, if $c\colonequals x+R\nu(x)$, then $B_R(c)\subset\R^3\setminus\overline\Omega$. These two balls are disjoint, and therefore
		$$
		|z-c|\geq R+\frac r2.
		$$
		It follows that
		$$
		|x-z|^2+2R(x-z)\cdot\nu(x)+R^2\geq R^2+Rr+\frac{r^2}{4}.
		$$
		Using \eqref{eq:geom-a} and $r\leq R/18$, we obtain
		$$
		(x-z)\cdot\nu(x)\geq\frac r2-\frac{|x-z|^2}{2R}\geq\frac r2-\frac{9r^2}{2R}\geq\frac r4.
		$$
	\end{enumerate}
\end{proof}

The following Pohozaev identity connects the energy in the interior of a ball with the energy on the spherical and physical parts of its boundary.

\begin{lemma}[Pohozaev identity up to the boundary]\label{lem:pohozaev-borde}
	Let $\Omega\subset\R^3$ be a bounded domain with $\partial\Omega\in C^2$ and outward unit normal $\nu$. Let $\ep>0$, $(v,\A)\in H^1(\Omega,\C)\times W^{1,6}(\Omega,\R^3)$ satisfy \eqref{eq:post_splitting_eq}, and assume that $\A$ is in the weighted Coulomb gauge \eqref{weightedCoulombgauge}. Set
	$$
	b_\ep\colonequals1-\frac{\ep^2|\A|^2}{\pmin^2},\qquad U\colonequals\Omega\cap B_r(x_0).
	$$
	Let $x_0\in\overline\Omega$, $z\in\R^3$, and $\eta(x)\colonequals x-z$. Denote by $\nu_r$ the outward unit normal to $B_r(x_0)$. Then, for almost every $r\in(0,\operatorname{diam}\Omega)$,
	\begin{align}\label{eq:poho-borde}
		&\frac12\int_U|\nabla v|^2+\frac{3}{4\ep^2}\int_U\pmin^2(b_\ep-|v|^2)^2\nonumber\\
		&\quad=\frac12\int_U\pmin^2(b_\ep-|v|^2)\eta\cdot\nabla\left(\frac{|\A|^2}{\pmin^2}\right)+2\int_U\ip{i\A\cdot\nabla v}{\eta\cdot\nabla v}\nonumber\\
		&\qquad-\int_U\frac2{\pmin}\ip{\nabla\pmin\cdot\nabla v}{\eta\cdot\nabla v}-\frac1{2\ep^2}\int_U\pmin(\eta\cdot\nabla\pmin)(b_\ep-|v|^2)^2\nonumber\\
		&\qquad+\int_{\Omega\cap\partial B_r}(\eta\cdot\nu_r)e_\ep-\int_{\Omega\cap\partial B_r}\left\langle\frac{\partial v}{\partial\nu_r},\eta\cdot\nabla v\right\rangle+\int_{\partial\Omega\cap B_r}(\eta\cdot\nu)e_\ep,
	\end{align}
	where
	$$
	e_\ep\colonequals\frac12|\nabla v|^2+\frac{\pmin^2}{4\ep^2}(b_\ep-|v|^2)^2
	$$
	and $i\A\cdot\nabla v\colonequals i\sum_k\A_k\partial_kv$.
\end{lemma}

\begin{remark}
	If $z=x_0$ and $B_r(x_0)\subset\Omega$, then $\eta=x-x_0$, $\eta\cdot\nu_r=r$, and $\eta\cdot\nabla v=r\partial_{\nu_r}v$ on $\partial B_r(x_0)$. The last integral in \eqref{eq:poho-borde} is absent, and the integrand on $\partial B_r(x_0)$ becomes
	$$
	\frac r2|\nabla_\top v|^2-\frac r2|\partial_{\nu_r}v|^2+\frac r{4\ep^2}\pmin^2(b_\ep-|v|^2)^2.
	$$
	This recovers the interior Pohozaev identity used, for example, in \cite{LiuZhou}. Thus, the case of balls contained in $\Omega$ is a simpler particular case of the identity above.
\end{remark}

\begin{proof}
	\begin{enumerate}[label=\textnormal{\bfseries Step~\arabic*.},ref=\arabic*,leftmargin=0pt,labelsep=*,itemindent=*,itemsep=10pt,topsep=10pt]
		\item \textbf{Flat form of the equation and regularity.} Using $\nabla_\A v=\nabla v-i\A v$, we have
		$$
		(\diver-i\A\cdot)(\pmin^2\nabla_\A v)=\diver(\pmin^2\nabla v)-i\diver(\pmin^2\A v)-i\pmin^2\A\cdot\nabla v-\pmin^2|\A|^2v.
		$$
		Since $\diver(\pmin^2\A)=0$,
		$$
		\diver(\pmin^2\A v)=\pmin^2\A\cdot\nabla v.
		$$
		Therefore, \eqref{eq:post_splitting_eq} can be written as
		\begin{equation}\label{eq:flat-v}
			\diver(\pmin^2\nabla v)=2i\pmin^2\A\cdot\nabla v+\pmin^4\partial_zW(x,v),
		\end{equation}
		where
		$$
		W(x,z)\colonequals\frac1{4\ep^2}(b_\ep(x)-|z|^2)^2,\qquad \partial_zW(x,z)=-\frac1{\ep^2}(b_\ep(x)-|z|^2)z.
		$$
		Indeed,
		$$
		\pmin^2|\A|^2v-\frac{\pmin^4}{\ep^2}(1-|v|^2)v=-\frac{\pmin^4}{\ep^2}\left(1-\frac{\ep^2|\A|^2}{\pmin^2}-|v|^2\right)v.
		$$
		We also note that $v\in H^2(\Omega)$, which justifies the integrations by parts below. Indeed, since $v\in H^1(\Omega)\hookrightarrow L^6(\Omega)$, we have $v(1-|v|^2)\in L^2(\Omega)$. Moreover, $\A\in W^{1,6}(\Omega,\R^3)\hookrightarrow L^\infty(\Omega,\R^3)$ and $\nabla\pmin\in L^\infty(\Omega,\R^3)$. Expanding the left-hand side of \eqref{eq:flat-v}, we may write
		$$
		\Delta v=\pmin^{-2}\left(2i\pmin^2\A\cdot\nabla v+\pmin^4\partial_zW(x,v)-\nabla\pmin^2\cdot\nabla v\right)\in L^2(\Omega).
		$$
		Moreover, $\partial_\nu v=0$ on $\partial\Omega$, since $\nabla_\A v\cdot\nu=0$ and $\A\cdot\nu=0$. Applying Theorem~\ref{thm:H2p-neumann} with $p=2$ to the real and imaginary parts gives $v\in H^2(\Omega)$.
		\item \textbf{Testing the equation.} Set
		$$
		h\colonequals\eta\cdot\nabla v.
		$$
		We multiply \eqref{eq:flat-v} by $h/\pmin^2$, take the real scalar product, and integrate over $U$:
		\begin{equation}\label{eq:poho-split}
			\underbrace{\int_U\ip{\diver(\pmin^2\nabla v)}{h/\pmin^2}}_{\mathrm I}=\underbrace{2\int_U\ip{i\A\cdot\nabla v}{h}}_{\mathrm{III}}+\underbrace{\int_U\pmin^2\ip{\partial_zW(x,v)}{h}}_{\mathrm{II}}.
		\end{equation}
		We use below that $\diver\eta=3$ and $\partial_k\eta_l=\delta_{kl}$.
		\item \textbf{The gradient term $\mathrm I$.} For $\xi,\zeta\in\C^3$ we denote by $\ip{\xi}{\zeta}\in\R$ the real inner product obtained by identifying $\C^3\cong\R^6$, that is, the one that extends the scalar product $\ip{\cdot}{\cdot}$ of $\C\cong\R^2$ defined in the introduction componentwise:
		\[
  			\ip{\xi}{\zeta}\coloneqq\sum_{i=1}^3\ip{\xi_i}{\zeta_i}=\operatorname{Re}\sum_{i=1}^3\overline{\xi_i}\,\zeta_i .
		\]
		In particular, $\ip{\xi}{\xi}=|\xi|^2$.
		Since $\partial_\nu v=0$ on $\partial\Omega$, integration by parts gives
		$$
		\mathrm I=-\int_U\pmin^2\ip{\nabla v}{\nabla(h/\pmin^2)}+\int_{\Omega\cap\partial B_r}\left\langle\frac{\partial v}{\partial\nu_r},h\right\rangle.
		$$
		By the quotient rule,
		$$
		\pmin^2\ip{\nabla v}{\nabla(h/\pmin^2)}=\ip{\nabla v}{\nabla h}-\frac2{\pmin}\ip{\nabla\pmin\cdot\nabla v}{\eta\cdot\nabla v}.
		$$
		Since $\partial_kh=\partial_kv+\eta\cdot\nabla(\partial_kv)$,
		$$
		\ip{\nabla v}{\nabla h}=|\nabla v|^2+\eta\cdot\nabla\frac{|\nabla v|^2}{2}.
		$$
		Moreover,
		$$
		\int_U\eta\cdot\nabla\frac{|\nabla v|^2}{2}=-\frac32\int_U|\nabla v|^2+\frac12\int_{\partial U}(\eta\cdot\nu_U)|\nabla v|^2.
		$$
		Consequently,
		\begin{multline}\label{eq:poho-I}
			\mathrm I=\frac12\int_U|\nabla v|^2+\int_U\frac2{\pmin}\ip{\nabla\pmin\cdot\nabla v}{\eta\cdot\nabla v}\\
			-\frac12\int_{\partial U}(\eta\cdot\nu_U)|\nabla v|^2+\int_{\Omega\cap\partial B_r}\left\langle\frac{\partial v}{\partial\nu_r},h\right\rangle.
		\end{multline}
		\item \textbf{The potential term $\mathrm{II}$.} Let $\mathcal W(x)\colonequals W(x,v(x))$. By the chain rule,
		$$
		\ip{\partial_zW(x,v)}{\eta\cdot\nabla v}=\eta\cdot\nabla\mathcal W-\eta\cdot\partial_xW(x,v).
		$$
		Since
		$$
		\diver(\pmin^2\eta)=3\pmin^2+2\pmin\eta\cdot\nabla\pmin,
		$$
		integration by parts yields
		\begin{align*}
			\int_U\pmin^2\eta\cdot\nabla\mathcal W&=-\frac3{4\ep^2}\int_U\pmin^2(b_\ep-|v|^2)^2-\frac1{2\ep^2}\int_U\pmin(\eta\cdot\nabla\pmin)(b_\ep-|v|^2)^2\\
			&\quad+\frac1{4\ep^2}\int_{\partial U}(\eta\cdot\nu_U)\pmin^2(b_\ep-|v|^2)^2.
		\end{align*}
		Furthermore,
		$$
		\partial_xW(x,v)=\frac{b_\ep-|v|^2}{2\ep^2}\nabla b_\ep,\qquad \nabla b_\ep=-\ep^2\nabla\left(\frac{|\A|^2}{\pmin^2}\right),
		$$
		and hence
		$$
		-\int_U\pmin^2\eta\cdot\partial_xW=\frac12\int_U\pmin^2(b_\ep-|v|^2)\eta\cdot\nabla\left(\frac{|\A|^2}{\pmin^2}\right).
		$$
		Therefore,
		\begin{align}\label{eq:poho-II}
			\mathrm{II}&=-\frac3{4\ep^2}\int_U\pmin^2(b_\ep-|v|^2)^2-\frac1{2\ep^2}\int_U\pmin(\eta\cdot\nabla\pmin)(b_\ep-|v|^2)^2\nonumber\\
			&\quad+\frac1{4\ep^2}\int_{\partial U}(\eta\cdot\nu_U)\pmin^2(b_\ep-|v|^2)^2+\frac12\int_U\pmin^2(b_\ep-|v|^2)\eta\cdot\nabla\left(\frac{|\A|^2}{\pmin^2}\right).
		\end{align}
		\item \textbf{Conclusion.} Substituting \eqref{eq:poho-I} and \eqref{eq:poho-II} into \eqref{eq:poho-split}, and rearranging the terms, we obtain
		\begin{align*}
			&\frac12\int_U|\nabla v|^2+\frac{3}{4\ep^2}\int_U\pmin^2(b_\ep-|v|^2)^2\\
			&\quad=\frac12\int_U\pmin^2(b_\ep-|v|^2)\eta\cdot\nabla\left(\frac{|\A|^2}{\pmin^2}\right)+2\int_U\ip{i\A\cdot\nabla v}{\eta\cdot\nabla v}\\
			&\qquad-\int_U\frac2{\pmin}\ip{\nabla\pmin\cdot\nabla v}{\eta\cdot\nabla v}-\frac1{2\ep^2}\int_U\pmin(\eta\cdot\nabla\pmin)(b_\ep-|v|^2)^2\\
			&\qquad+\int_{\partial U}(\eta\cdot\nu_U)e_\ep-\int_{\Omega\cap\partial B_r}\left\langle\frac{\partial v}{\partial\nu_r},\eta\cdot\nabla v\right\rangle.
		\end{align*}
		Splitting
		$$
		\partial U=(\Omega\cap\partial B_r)\cup(\partial\Omega\cap B_r)
		$$
		and recalling that $\nu_U=\nu_r$ on the first part and $\nu_U=\nu$ on the second gives \eqref{eq:poho-borde}.
	\end{enumerate}
\end{proof}

In the remainder of this section, we assume that there exist $M_A>0$ and $C_0>0$ such that
\begin{equation}\label{eq:trace-hyp}
	\|\A\|_{L^\infty(\Omega,\R^3)}+\|\A\|_{W^{1,6}(\Omega,\R^3)}\leq C_0\lep^{M_A}.
\end{equation}
Increasing $M_A$ if necessary, we assume without loss of generality that $M_A\geq M_\infty$. Notice that since $W^{1,6}(\Omega,\R^3)\hookrightarrow L^\infty(\Omega,\R^3)$, the first term on the left-hand side is controlled by the second one; we retain both for clarity.

We also assume that $\pmin$ and $(v,\A)$ satisfy the hypotheses of Lemma~\ref{lem:pohozaev-borde}, and in particular that $\pmin$ satisfies \eqref{regularitypmin}. For $x_0\in\overline\Omega$, we set
\begin{equation}\label{eq:definitions_monotonicity}
	E(r)\colonequals\int_{\Omega\cap B_r(x_0)}e_\ep,\qquad \Phi(r)\colonequals\frac{E(r)}r,\qquad K_\ep\colonequals\ep^2|\log\ep|^{5M_A}.
\end{equation}

The next lemma controls the four interior error terms arising on the right-hand side of the Pohozaev identity.

\begin{lemma}[Error terms in the Pohozaev identity]\label{lem:poho-err}
	Let $x_0\in\overline\Omega$, $0<r<1$, $U\colonequals\Omega\cap B_r(x_0)$, and let $z\in\R^3$ be such that $\eta(x)\colonequals x-z$ satisfies $|\eta|\leq3r$ in $U$. Let $\mathrm R$ denote the sum of the first four integrals on the right-hand side of \eqref{eq:poho-borde}. Then
	\begin{equation}\label{eq:poho-err}
		|\mathrm R|\leq C rE(r)\lep^{M_A}+Cr^2K_\ep.
	\end{equation}
\end{lemma}

\begin{proof}
	Since $b\leq\pmin^2\leq1$, we have $\pmin\leq\pmin^2/\sqrt b$. Moreover, by the definition of $e_\ep$,
	\begin{equation}\label{eq:err-basic}
		\int_U|\nabla v|^2\leq2E(r),\qquad \int_U\pmin^2(b_\ep-|v|^2)^2\leq4\ep^2E(r),\qquad |U|\leq\frac{4\pi}{3}r^3.
	\end{equation}
	\begin{enumerate}[label=\textnormal{\bfseries Step~\arabic*.},ref=\arabic*,leftmargin=0pt,labelsep=*,itemindent=*,itemsep=10pt,topsep=10pt]
		\item\label{step1errorpoh} \textbf{Terms involving $\A$, $\nabla\pmin$, and $\nabla v$.} By \eqref{eq:trace-hyp} and \eqref{eq:err-basic},
		\begin{align*}
			\left|2\int_U\ip{i\A\cdot\nabla v}{\eta\cdot\nabla v}\right|&\leq6r\|\A\|_{L^\infty}\int_U|\nabla v|^2\leq C rE(r)\lep^{M_A},\\
			\left|\int_U\frac2{\pmin}\ip{\nabla\pmin\cdot\nabla v}{\eta\cdot\nabla v}\right|&\leq\frac{6r}{\sqrt b}\|\nabla\pmin\|_{L^\infty}\int_U|\nabla v|^2\leq C rE(r)\lep^{M_A},\\
			\frac1{2\ep^2}\left|\int_U\pmin(\eta\cdot\nabla\pmin)(b_\ep-|v|^2)^2\right|&\leq\frac{3r\|\nabla\pmin\|_{L^\infty}}{2\ep^2\sqrt b}\int_U\pmin^2(b_\ep-|v|^2)^2\\
			&\leq CrE(r)\lep^{M_A}.
		\end{align*}
		\item \textbf{The term involving $\nabla(|\A|^2/\pmin^2)$.} Since
		$$
		\nabla\left(\frac{|\A|^2}{\pmin^2}\right)=\frac2{\pmin^2}\sum_j\A_j\nabla\A_j-\frac{2|\A|^2}{\pmin^3}\nabla\pmin,
		$$
		\eqref{eq:trace-hyp} gives
		$$
		\left|\nabla\left(\frac{|\A|^2}{\pmin^2}\right)\right|\leq C\lep^{M_A}|\nabla\A|+C\lep^{3M_A}.
		$$
		Thus, the first term in $\mathrm R$ is bounded by
		\begin{equation}\label{eq:bound_first_R}
			Cr\int_U\pmin|b_\ep-|v|^2|\left(\lep^{M_A}|\nabla\A|+\lep^{3M_A}\right).
		\end{equation}
		By the Cauchy--Schwarz and H\"older inequalities,
		$$
		\|\nabla\A\|_{L^2(U)}\leq|U|^{1/3}\|\nabla\A\|_{L^6(\Omega)}\leq Cr\lep^{M_A},
		$$
		and hence
		\begin{align*}
			\int_U\pmin|b_\ep-|v|^2||\nabla\A|&\leq2C\ep r\sqrt{E(r)}\,\lep^{M_A},\\
			\int_U\pmin|b_\ep-|v|^2|&\leq2\ep\sqrt{E(r)}\,|U|^{1/2}\leq C\ep r^{3/2}\sqrt{E(r)}.
		\end{align*}
		Since $r<1$ and $\lep\geq1$, \eqref{eq:bound_first_R} is bounded by
		$$
		Cr^2\ep\lep^{3M_A}\sqrt{E(r)}.
		$$
		Young's inequality then gives
		\begin{align*}
			r^2\ep\lep^{3M_A}\sqrt{E(r)}&=\left(\sqrt{\lep^{M_A} rE(r)}\right)\big(r^{3/2}\ep\lep^{\tfrac52M_A}\big)\\
			&\leq\frac12\lep^{M_A} rE(r)+\frac12r^3\ep^2\lep^{5M_A}\\
			&\leq\frac12\lep^{M_A} rE(r)+\frac12r^2K_\ep.
		\end{align*}
		Combining this estimate with Step~\ref{step1errorpoh} proves \eqref{eq:poho-err}.
	\end{enumerate}
\end{proof}

The next lemma uses the Pohozaev identity with a suitably shifted center. The geometric choice of this center makes the contribution on $\partial\Omega$ positive, allowing us to control the trace of the energy density by interior and spherical terms.

\begin{lemma}[Trace estimate]\label{lem:trace}
	Let $r_1$ be as in Lemma~\ref{lem:geom-borde}. There exists $C>0$ such that, for every $x_0\in\overline\Omega$, setting $d_0\colonequals\dist(x_0,\partial\Omega)$, and for almost every $r\in[d_0,r_1)$,
	\begin{equation}\label{eq:trace}
		\int_{\partial\Omega\cap B_r(x_0)}e_\ep\leq C\left(1+r\lep^{M_A}\right)\Phi(r)+CE'(r)+CrK_\ep.
	\end{equation}
\end{lemma}

\begin{proof}
	\begin{enumerate}[label=\textnormal{\bfseries Step~\arabic*.},ref=\arabic*,leftmargin=0pt,labelsep=*,itemindent=*,itemsep=10pt,topsep=10pt]
		\item \textbf{Choice of the center.} Let $p\in\partial\Omega$ satisfy $|x_0-p|=d_0$, and set
		$$
		z\colonequals p-\frac r2\nu(p),\qquad \eta(x)\colonequals x-z,\qquad U\colonequals\Omega\cap B_r(x_0).
		$$
		By \eqref{eq:geom-a} and \eqref{eq:geom-b},
		\begin{equation}\label{eq:geometry_for_trace}
			|\eta|\leq3r\quad\text{in }\overline{B_r(x_0)},\qquad \eta\cdot\nu\geq\frac r4\quad\text{on }\partial\Omega\cap B_r(x_0).
		\end{equation}
		Writing \eqref{eq:poho-borde} as
		\begin{equation*}
			\int_{\partial\Omega\cap B_r(x_0)}(\eta\cdot\nu)e_\ep=\mathrm L-\mathrm R-\mathrm S,
		\end{equation*}
		where $\mathrm L$ denotes the left-hand side of \eqref{eq:poho-borde}, $\mathrm R$ is as in Lemma~\ref{lem:poho-err}, and
		$$
		\mathrm S\colonequals\int_{\Omega\cap\partial B_r(x_0)}\left[(\eta\cdot\nu_r)e_\ep-\left\langle\frac{\partial v}{\partial\nu_r},\eta\cdot\nabla v\right\rangle\right],
		$$
		we obtain from \eqref{eq:geometry_for_trace}
		\begin{equation}\label{eq:trace-start}
			\frac r4\int_{\partial\Omega\cap B_r(x_0)}e_\ep\leq\mathrm L+|\mathrm R|+|\mathrm S|.
		\end{equation}
		\item \textbf{Interior and spherical contributions.} From the definition of $e_\ep$,
		$$
		\mathrm L=\frac12\int_U|\nabla v|^2+\frac{3}{4\ep^2}\int_U\pmin^2(b_\ep-|v|^2)^2\leq3E(r).
		$$
		On $\Omega\cap\partial B_r$, \eqref{eq:geometry_for_trace} gives $|\eta\cdot\nu_r|\leq3r$, while
		$$
		|\partial_{\nu_r}v|\,|\eta\cdot\nabla v|\leq3r|\nabla v|^2\leq6re_\ep.
		$$
		Therefore,
		$$
		|\mathrm S|\leq9r\int_{\Omega\cap\partial B_r(x_0)}e_\ep=9rE'(r)
		$$
		for almost every $r$.
		\item \textbf{Conclusion.} Lemma~\ref{lem:poho-err} and \eqref{eq:trace-start} yield
		$$
		\frac r4\int_{\partial\Omega\cap B_r(x_0)}e_\ep\leq3E(r)+C\lep^{M_A} rE(r)+Cr^2K_\ep+9rE'(r).
		$$
		Dividing by $r/4$ and using $E(r)=r\Phi(r)$ proves \eqref{eq:trace}.
	\end{enumerate}
\end{proof}

We now combine the Pohozaev identity, the interior error estimates, and the trace estimate to transport the normalized energy between different radii.

\begin{lemma}[Pseudomonotonicity up to the boundary]\label{lem:mono-borde}
	Let $r_1$ and $\kappa$ be as in Lemma~\ref{lem:geom-borde}, and let $M_A$, $K_\ep$, $E$, and $\Phi$ be as above, so that \eqref{eq:trace-hyp} holds. There exists $C>0$ such that, for every $x_0\in\overline\Omega$ and every $0<s<r<r_1$,
	\begin{equation}\label{eq:mono}
		\Phi(s)+K_\ep\leq e^{C(r\lep^{M_A} +1)}\big(\Phi(r)+K_\ep\big).
	\end{equation}
\end{lemma}

\begin{proof}
	Fix $x_0\in\overline\Omega$, set $d_0\colonequals\dist(x_0,\partial\Omega)$, and write $B_r\colonequals B_r(x_0)$. The function $E$ is absolutely continuous, and, for almost every $r\in(0,r_1)$,
	\begin{equation}\label{eq:mono-E}
		E'(r)=\int_{\Omega\cap\partial B_r}e_\ep=\Phi(r)+r\Phi'(r),\qquad \Phi'(r)=\frac{rE'(r)-E(r)}{r^2}.
	\end{equation}
	\begin{enumerate}[label=\textnormal{\bfseries Step~\arabic*.},ref=\arabic*,leftmargin=0pt,labelsep=*,itemindent=*,itemsep=10pt,topsep=10pt]
		\item \textbf{Differential inequality from the Pohozaev identity.} Apply Lemma~\ref{lem:pohozaev-borde} with $z=x_0$, so that $\eta(x)=x-x_0$. On $\Omega\cap\partial B_r$, we have $\eta=r\nu_r$ and hence
		\begin{align}\label{eq:mono_bound_E'_ball}
			&\int_{\Omega\cap\partial B_r}(\eta\cdot\nu_r)e_\ep-\int_{\Omega\cap\partial B_r}\left\langle\frac{\partial v}{\partial\nu_r},\eta\cdot\nabla v\right\rangle\nonumber\\
			&\qquad=rE'(r)-r\int_{\Omega\cap\partial B_r}|\partial_{\nu_r}v|^2\leq rE'(r).
		\end{align}
		Moreover,
		\begin{equation}\label{eq:mono_E_poho}
			E(r)\leq\frac12\int_{\Omega\cap B_r}|\nabla v|^2+\frac{3}{4\ep^2}\int_{\Omega\cap B_r}\pmin^2(b_\ep-|v|^2)^2.
		\end{equation}
		Using \eqref{eq:mono-E}, \eqref{eq:mono_bound_E'_ball}, \eqref{eq:mono_E_poho}, the Pohozaev identity, and Lemma~\ref{lem:poho-err}, we obtain
		\begin{equation}\label{eq:mono_E_pho_complete_chad}
			\Phi'(r)\geq-C\lep^{M_A}\Phi(r)-CK_\ep-\frac1{r^2}\int_{\partial\Omega\cap B_r}|(x-x_0)\cdot\nu|e_\ep.
		\end{equation}
		This is the basic differential inequality. It remains to control its boundary term.
		\item \textbf{Radii below the distance to the boundary.} If $r<d_0$, then $B_r(x_0)\subset\Omega$, and the boundary term in \eqref{eq:mono_E_pho_complete_chad} is absent. Since $\lep\geq1$,
		\begin{equation}\label{eq:mono-int}
			(\Phi+K_\ep)'(r)=\Phi'(r)\geq-C\lep^{M_A}\Phi(r)-CK_\ep\geq-C\lep^{M_A}(\Phi(r)+K_\ep).
		\end{equation}
		\item \textbf{Radii intersecting the boundary.} Suppose that $d_0\leq r<r_1$. By \eqref{eq:geom-b},
		$$
		|(x-x_0)\cdot\nu|\leq d_0+4\kappa r^2=ar^2\quad\text{on }\partial\Omega\cap B_r,
		$$
		where
		$$
		a\colonequals\frac{d_0}{r^2}+4\kappa.
		$$
		Notice that $ar\leq1+4\kappa$. Lemma~\ref{lem:trace} and \eqref{eq:mono-E} give
		\begin{align*}
			\frac1{r^2}\int_{\partial\Omega\cap B_r}|(x-x_0)\cdot\nu|e_\ep&\leq a\int_{\partial\Omega\cap B_r}e_\ep\nonumber\\
			&\leq a\left[C(1+\lep^{M_A} r)\Phi+C(\Phi+r\Phi')+CrK_\ep\right].
		\end{align*}
		Substituting this estimate into \eqref{eq:mono_E_pho_complete_chad} and using $a\lep^{M_A} r\leq C\lep^{M_A}$, $ar\leq C$, and $\lep\geq1$, we obtain
		$$
		(1+Car)\Phi'(r)\geq-C(\lep^{M_A}+a)\Phi(r)-CK_\ep.
		$$
		Since $1+Car$ is bounded above and below by positive constants and
		$$
		a=\frac{d_0}{r^2}+4\kappa,
		$$
		we conclude, after absorbing $4\kappa$ into $\lep^{M_A}\geq1$, that
		\begin{equation}\label{eq:mono-bd}
			(\Phi+K_\ep)'(r)\geq-C\left(\lep^{M_A}+\frac{d_0}{r^2}\right)(\Phi(r)+K_\ep)
		\end{equation}
		for almost every $r\in[d_0,r_1)$.
		\item \textbf{Integration of the differential inequality.} Combining \eqref{eq:mono-int} and \eqref{eq:mono-bd}, we have
		$$
		(\Phi+K_\ep)'(t)\geq-c(t)(\Phi(t)+K_\ep)
		$$
		for almost every $t\in(0,r_1)$, where
		$$
		c(t)\colonequals C\lep^{M_A}+C\frac{d_0}{t^2}\chi_{\{t\geq d_0\}}.
		$$
		Thus,
		$$
		t\longmapsto\exp\left(\int_0^tc(\tau)\,d\tau\right)(\Phi(t)+K_\ep)
		$$
		is nondecreasing. Hence, for $0<s<r<r_1$,
		$$
		(\Phi+K_\ep)(s)\leq\exp\left(\int_s^rc(t)\,dt\right)(\Phi+K_\ep)(r).
		$$
		Moreover,
		$$
		\int_s^rc(t)\,dt\leq Cr\lep^{M_A}+C\int_{\max\{s,d_0\}}^\infty\frac{d_0}{t^2}\,dt\leq C(r\lep^{M_A} +1),
		$$
		because $d_0/\max\{s,d_0\}\leq1$. This proves \eqref{eq:mono}.
	\end{enumerate}
\end{proof}

We conclude with a clearing-out estimate adapted to the pinned energy density $e_\ep$. This is the small-scale ingredient that will later be combined with Lemma~\ref{lem:mono-borde} to transfer energy control from larger scales down to balls of radius $\ep$.

\begin{lemma}[Pinned clearing-out estimate]\label{lem:eps-reg}
	Let $M_A$, $E$, and $\Phi$ be as above, so that \eqref{eq:trace-hyp} holds. Assume in addition that $|v|\leq1$ in $\Omega$ and that, for some $C_1>0$,
	\begin{equation}\label{eq:grad-hyp}
		\|\nabla v\|_{L^\infty(\Omega,\C^3)}\leq\frac{C_1}{\ep}.
	\end{equation}
	Then there exists $C>0$, depending only on $\Omega$, $b$, $C_0$, and $C_1$, such that, for every $x_0\in\overline\Omega$,
	\begin{equation}\label{eq:eps-reg}
		(1-|v(x_0)|)^5\leq C\big(\Phi(\ep)+\ep^4\lep^{4M_A}\big),\qquad \Phi(\ep)=\frac1\ep\int_{\Omega\cap B_\ep(x_0)}e_\ep.
	\end{equation}
\end{lemma}

\begin{proof}
	Since $|\nabla|v||\leq|\nabla v|$, assumption \eqref{eq:grad-hyp} allows us to apply Lemma~\ref{lemma:clearingout} to $v$, yielding
	\begin{equation}\label{eq:eps-reg-1}
		(1-|v(x_0)|)^5\leq\frac C{\ep^3}\int_{\Omega\cap B_\ep(x_0)}(1-|v|^2)^2.
	\end{equation}
	It remains to compare $(1-|v|^2)^2$ with the potential term in $e_\ep$. We write
	$$
	1-|v|^2=(b_\ep-|v|^2)+(1-b_\ep).
	$$
	By \eqref{eq:trace-hyp} and $\pmin^2\geq b$,
	$$
	0\leq1-b_\ep=\frac{\ep^2|\A|^2}{\pmin^2}\leq\frac{C_0^2}{b}\ep^2\lep^{2M_A}.
	$$
	Therefore, using $(\alpha+\beta)^2\leq2\alpha^2+2\beta^2$, $\pmin^2\geq b$, and for $U\colonequals\Omega\cap B_\ep(x_0)$, $|U|\leq |B_\ep(x_0)|=\tfrac{4\pi\ep^3}{3}$, we obtain
	\begin{align*}
		\int_U(1-|v|^2)^2&\leq\frac2b\int_U\pmin^2(b_\ep-|v|^2)^2+2|U|\frac{C_0^4}{b^2}\ep^4\lep^{4M_A}\\
		&\leq\frac{8\ep^2}{b}E(\ep)+C\ep^7\lep^{4M_A}.
	\end{align*}
	Multiplying by $\ep^{-3}$ and using $E(\ep)=\ep\Phi(\ep)$ gives
	$$
	\frac1{\ep^3}\int_U(1-|v|^2)^2\leq C\Phi(\ep)+C\ep^4\lep^{4M_A}.
	$$
	Combining this estimate with \eqref{eq:eps-reg-1} proves \eqref{eq:eps-reg}.
\end{proof}

\section{Proof of Theorem~\ref{thm:main}}\label{sec:proofthm1}
This section is devoted to the proof of Theorem~\ref{thm:main}. The argument follows closely the proof of \cite[Theorem 1.2]{diaz-roman-3d}, with the main additional ingredient being the logarithmic-loss H\"older estimate for $\Bfield$ established in Proposition~\ref{prop:Bfield-holder}.

\begin{proof}[Proof of Theorem~\ref{thm:main}]
	We divide the proof into three steps. Step~1 establishes an a priori bound for the weighted free energy, Step~2 proves item~\ref{thm1:itema}, and Step~3 proves items~\ref{thm1:itemb} and \ref{thm1:itemc}.
	\begin{enumerate}[label=\textnormal{\bfseries Step~\arabic*.},ref=\arabic*,leftmargin=0pt,labelsep=*,itemindent=*,itemsep=10pt,topsep=10pt]
		\item \textbf{A priori bound for the weighted free energy.} Since $(\u,\A)$ minimizes $\glthreed$, we have
		$$
		\glthreed(\u,\A)\leq\glthreed\big(\pmin e^{ih_\ex\phi_\ep^0},h_\ex A_\ep^0\big).
		$$
		Thus, by the energy splitting \eqref{Split3D},
		\begin{equation}\label{eq:upper-split}
			\fenweightthreed(u,A)\leq h_\ex\int_\Omega\mu(u,A)\cdot\Bfield+|R_0|.
		\end{equation}
		By integration by parts,
		$$
		\int_\Omega\mu(u,A)\cdot\Bfield=\int_\Omega\ip{iu}{\nabla_Au}\cdot\curl\Bfield+\int_\Omega\curl A\cdot\Bfield.
		$$
		Hence
		\begin{align}
			\fenweightthreed(u,A)&\leq h_\ex\int_\Omega\left(\ip{iu}{\nabla_Au}\cdot\curl\Bfield+\curl A\cdot\Bfield\right)+|R_0|\notag\\
			&\leq C\left(h_\ex\|u\|_{L^\infty}\|\nabla_Au\|_{L^2}\|\curl\Bfield\|_{L^2}+h_\ex\|\curl A\|_{L^2}\|\Bfield\|_{L^2}+|R_0|\right)\notag\\
			&\leq C\left(h_\ex\|\nabla_Au\|_{L^2}+h_\ex\|\curl A\|_{L^2}+|R_0|\right),\label{eq:step1-holder}
		\end{align}
		where we used \eqref{boundu}, \eqref{estimateforA0}, and the Friedrichs--Poincar\'e inequality. Moreover,
		\begin{equation}\label{eq:grad-curl}
			\|\nabla_Au\|_{L^2(\Omega,\C^3)}+\|\curl A\|_{L^2(\Omega,\R^3)}\leq C\fenweightthreed(u,A)^{1/2}.
		\end{equation}
		Using $0\leq1-|u|^2\leq1$, $\pmin^2\geq b$, \eqref{estimateforcurlB0Linf}, and the Cauchy--Schwarz inequality, we obtain
		\begin{align}
			|R_0|&\leq\frac{h_\ex^2}{2}\int_\Omega\frac{|\curl\Bfield|^2}{\pmin^2}|1-|u|^2|\notag\\
			&\leq Ch_\ex^2\|\curl\Bfield\|_{L^\infty}^2\left(\int_\Omega(1-|u|^2)^2\right)^{1/2}\notag\\
			&\leq Ch_\ex^2|\log\ep|^{2\McurlBfield}\ep\,\fenweightthreed(u,A)^{1/2},\label{eq:r-bound}
		\end{align}
		where we used
		$$
		\int_\Omega(1-|u|^2)^2\leq C\ep^2\fenweightthreed(u,A).
		$$
		Since $h_\ex=O(|\log\ep|)$,
		$$
		|R_0|\leq Ch_\ex\big(\ep|\log\ep|^{1+2\McurlBfield}\big)\fenweightthreed(u,A)^{1/2}\leq Ch_\ex\fenweightthreed(u,A)^{1/2}
		$$
		for $\ep$ sufficiently small. Combining \eqref{eq:step1-holder}, \eqref{eq:grad-curl}, and \eqref{eq:r-bound}, we conclude that
		$$
		\fenweightthreed(u,A)=O\left(h_\ex\fenweightthreed(u,A)^{1/2}\right),
		$$
		and therefore
		\begin{equation}\label{eq:F-bound-hex}
			\fenweightthreed(u,A)\leq Ch_\ex^2.
		\end{equation}
		\item \textbf{Vorticity estimate at exponent $\gamma_0$.} Since $h_\ex=O(|\log\ep|)$, \eqref{eq:F-bound-hex} yields
		\begin{equation}\label{eq:F-log2}
			\fenweightthreed(u,A)\leq C|\log\ep|^2.
		\end{equation}
		Therefore, Theorem~\ref{epleveltools} applies with $m=2$. Let $\gamma_0$ be as in the statement and choose $n\in\N$ satisfying
		\begin{equation}\label{eq:choice-n}
			n>N+\McurlBfield+2\qquad\text{and}\qquad n>\frac{N+2\McurlBfield+5}{3\gamma_0}-2.
		\end{equation}
		Let $q\colonequals q(2,n)=3(n+2)$ and let $\nu_\ep$ be the current furnished by Theorem~\ref{epleveltools}. Using \eqref{vorticity estimate} and \eqref{estimateforB0holder}, we obtain
		\begin{align}
			\int_\Omega\mu(u,A)\cdot\Bfield&\leq\|\mu(u,A)-\nu_\ep\|_{(C_T^{0,\gamma_0}(\Omega,\R^3))^*}\|\Bfield\|_{C_T^{0,\gamma_0}(\Omega,\R^3)}+\langle\nu_\ep,\Bfield\rangle\notag\\
			&\leq\langle\nu_\ep,\Bfield\rangle+C|\log\ep|^{2+\McurlBfield-\gamma_0q}.\label{eq:mu-vs-nu}
		\end{align}
		Multiplying \eqref{eq:mu-vs-nu} by $-h_\ex$ and using $h_\ex\leq C|\log\ep|$, we get
		$$
		-h_\ex\int_\Omega\mu(u,A)\cdot\Bfield\geq-h_\ex\langle\nu_\ep,\Bfield\rangle-C|\log\ep|^{3+\McurlBfield-\gamma_0q}.
		$$
		On the other hand, \eqref{weighted lower bound} and the definition of $\fenweightthreed$ give
		\begin{align*}
			\fenweightthreed(u,A)&\geq\frac12\int_{S_{\nu_\ep}}\left(\pmin^2|\nabla_Au|^2+\frac{\pmin^4}{2\ep^2}(1-|u|^2)^2+|\curl A|^2\right)\\
			&\geq\frac{|\pmin^2\nu_\ep|}{2}\big(|\log\ep|-C\log|\log\ep|\big)-\frac{C}{|\log\ep|^n}.
		\end{align*}
		Combining the preceding two estimates, we find that
		\begin{align*}
			\fenweightthreed(u,A)-h_\ex\int_\Omega\mu(u,A)\cdot\Bfield
			&\geq\frac{|\pmin^2\nu_\ep|}{2}\big(|\log\ep|-C\log|\log\ep|\big)-h_\ex\langle\nu_\ep,\Bfield\rangle\\
			&\quad-\frac{C}{|\log\ep|^n}-C|\log\ep|^{3+\McurlBfield-\gamma_0q}.
		\end{align*}
		By \eqref{eq:choice-n},
		$$
		\frac1{|\log\ep|^n}+|\log\ep|^{3+\McurlBfield-\gamma_0q}=o\big(|\log\ep|^{-N-\McurlBfield-2}\big),
		$$
		and therefore
		\begin{equation}\label{eq:lb-combined}
			\begin{aligned}
				\fenweightthreed(u,A)-h_\ex\int_\Omega\mu(u,A)\cdot\Bfield
				&\geq\frac{|\pmin^2\nu_\ep|}{2}\big(|\log\ep|-C\log|\log\ep|\big)-h_\ex\langle\nu_\ep,\Bfield\rangle\\
				&\quad+o(|\log\ep|^{-N-\McurlBfield-2}).
			\end{aligned}
		\end{equation}
		By the definition of $\mathrm R_\ep$,
		$$
		\langle\nu_\ep,\Bfield\rangle\leq|\pmin^2\nu_\ep|\mathrm R_\ep.
		$$
		Inserting this estimate into \eqref{eq:lb-combined}, we obtain
		\begin{align*}
			\fenweightthreed(u,A)-h_\ex\int_\Omega\mu(u,A)\cdot\Bfield
			&\geq\frac{|\pmin^2\nu_\ep|}{2}\big(|\log\ep|-2h_\ex\mathrm R_\ep-C\log|\log\ep|\big)\\
			&\quad+o(|\log\ep|^{-N-\McurlBfield-2}).
		\end{align*}
		Since $\critfield=|\log\ep|/(2\mathrm R_\ep)$ and $h_\ex\leq\critfield-K_0\log|\log\ep|$, we have
		$$
		|\log\ep|-2h_\ex\mathrm R_\ep\geq2K_0\mathrm R_\ep\log|\log\ep|.
		$$
		Consequently,
		\begin{align*}
			\fenweightthreed(u,A)-h_\ex\int_\Omega\mu(u,A)\cdot\Bfield
			&\geq\frac{|\pmin^2\nu_\ep|}{2}(2K_0\mathrm R_\ep-C)\log|\log\ep|\\
			&\quad+o(|\log\ep|^{-N-\McurlBfield-2}).
		\end{align*}
		Since $\liminf_{\ep\to0^+}\mathrm R_\ep>0$, we may choose $K_0$, independently of $\ep$, such that $2K_0\mathrm R_\ep-C>1$ for $\ep$ sufficiently small. Hence
		\begin{equation}\label{eq:lb-loglog}
			\fenweightthreed(u,A)-h_\ex\int_\Omega\mu(u,A)\cdot\Bfield\geq\frac12|\pmin^2\nu_\ep|\log|\log\ep|+o(|\log\ep|^{-N-\McurlBfield-2}).
		\end{equation}
		On the other hand, \eqref{eq:upper-split}, \eqref{eq:r-bound}, and \eqref{eq:F-log2} imply
		\begin{equation}\label{eq:ub-r}
			\fenweightthreed(u,A)-h_\ex\int_\Omega\mu(u,A)\cdot\Bfield\leq|R_0|=O\left(\ep|\log\ep|^{3+2\McurlBfield}\right).
		\end{equation}
		Combining \eqref{eq:lb-loglog} and \eqref{eq:ub-r}, we deduce that
		$$
		|\pmin^2\nu_\ep|=o(|\log\ep|^{-N-\McurlBfield-1}).
		$$
		Since $b\leq\pmin^2\leq1$,
		$$
		|\nu_\ep|\leq\frac1b|\pmin^2\nu_\ep|=o(|\log\ep|^{-N-\McurlBfield-1}),
		$$
		and hence
		$$
		\|\nu_\ep\|_{(C_T^{0,\gamma_0}(\Omega,\R^3))^*}\leq|\nu_\ep|=o(|\log\ep|^{-N-\McurlBfield-1}).
		$$
		Using once more \eqref{vorticity estimate}, \eqref{eq:F-log2}, and \eqref{eq:choice-n}, we obtain
		\begin{equation}\label{eq:mu-small}
			\|\mu(u,A)\|_{(C_T^{0,\gamma_0}(\Omega,\R^3))^*}=o(|\log\ep|^{-N-\McurlBfield-1}).
		\end{equation}
		This proves item~\ref{thm1:itema} for $\gamma=\gamma_0$. If $\gamma\in[\gamma_0,1]$, then
		$$
		C_T^{0,\gamma}(\Omega,\R^3)\hookrightarrow C_T^{0,\gamma_0}(\Omega,\R^3),
		$$
		and therefore
		$$
		\|\mu(u,A)\|_{(C_T^{0,\gamma}(\Omega,\R^3))^*}\leq C\|\mu(u,A)\|_{(C_T^{0,\gamma_0}(\Omega,\R^3))^*}.
		$$
		Thus, item~\ref{thm1:itema} holds for every $\gamma\in[\gamma_0,1]$.
		\item \textbf{Approximation by the Meissner state.} Using \eqref{eq:mu-small} and \eqref{estimateforB0holder} with $\gamma=\gamma_0$, we obtain
		\begin{align*}
			h_\ex\left|\int_\Omega\mu(u,A)\cdot\Bfield\right|
			&\leq h_\ex\|\mu(u,A)\|_{(C_T^{0,\gamma_0}(\Omega,\R^3))^*}\|\Bfield\|_{C_T^{0,\gamma_0}(\Omega,\R^3)}\\
			&\leq C|\log\ep|\,o(|\log\ep|^{-N-\McurlBfield-1})\,|\log\ep|^{\McurlBfield}\\
			&=o(|\log\ep|^{-N}).
		\end{align*}
		Moreover, \eqref{eq:ub-r} implies $|R_0|=o(|\log\ep|^{-N})$ for every $N\in\N$. Therefore, \eqref{eq:upper-split} gives
		\begin{equation}\label{eq:F-small}
			\fenweightthreed(u,A)=o(|\log\ep|^{-N}),
		\end{equation}
		which proves item~\ref{thm1:itemb}. Finally, inserting \eqref{eq:F-small}, \eqref{eq:mu-small}, and \eqref{eq:r-bound} into the energy splitting \eqref{Split3D}, we conclude that
		$$
		\glthreed(\u,\A)=\glthreed\big(\pmin e^{ih_\ex\phi_\ep^0},h_\ex A_\ep^0\big)+o(|\log\ep|^{-N}),
		$$
		which proves item~\ref{thm1:itemc}.
	\end{enumerate}
	The theorem is thus proved.
\end{proof}

\section{Proof of Theorem~\ref{thm:clearing-final}}\label{sec:proofthm2}
This section is devoted to the proof of Theorem~\ref{thm:clearing-final}. The main task is to relate the small free energy of the fluctuation around the Meissner state to the local energy of a gauge-equivalent total configuration satisfying the equation required by the Pohozaev identity.

\begin{proof}[Proof of Theorem~\ref{thm:clearing-final}]
	Throughout the proof, $C>0$ denotes constants independent of $\ep$ and $x_0$, whose value may change from line to line. By Theorem~\ref{thm:main}, applied with $N=k+1$,
	\begin{equation}\label{eq:F-small-final}
		\fenweightthreed(u,A)=o\big(\lep^{-k-1}\big).
	\end{equation}
	Since $\fenweightthreed$ and $\curl A$ are gauge invariant, \eqref{eq:F-small-final} holds in any gauge. The pair $(u,A)$ is the fluctuation around the Meissner state, for which \eqref{eq:F-small-final} holds, whereas the equation \eqref{eq:post_splitting_eq} required by the Pohozaev identity is satisfied by the total pair. We therefore work with both.
	\begin{enumerate}[label=\textnormal{\bfseries Step~\arabic*.},ref=\arabic*,leftmargin=0pt,labelsep=*,itemindent=*,itemsep=10pt,topsep=10pt]
		\item\label{step:total-weighted-gauge} \textbf{The total pair in the weighted Coulomb gauge.} Since a global minimizer is a critical point, the pair $(\pmin^{-1}\u,\A)$ satisfies \eqref{eq:post_splitting_eq}. Moreover, Proposition~\ref{prop:regularity} gives $|\u|\leq1$ and $|\u|/\pmin\leq1$. Let $\varphiCouGauge$ be a phase putting the total pair in the Coulomb gauge in $\Omega$, and write
		$$
		\ACouGauge\colonequals\A-\nabla\varphiCouGauge,
		$$
		with $\uCouGauge$ denoting the corresponding reduced order parameter. By Proposition~\ref{prop:regularityA},
		$$
		\|\ACouGauge\|_{W^{1,6}(\Omega,\R^3)}\leq C\left(\|\nabla_{\ACouGauge}\uCouGauge\|_{L^2(\Omega,\C^3)}+h_\ex\right).
		$$
		The $L^2$ norm of the covariant gradient is gauge invariant. Since $\u/\pmin=u e^{ih_\ex\phi_\ep^0}$, a direct computation gives
		$$
		\left|\nabla_{A+h_\ex A_\ep^0}(\u/\pmin)\right|=\left|\nabla_Au-ih_\ex(A_\ep^0-\nabla\phi_\ep^0)u\right|.
		$$
		By \eqref{defBfield},
		$$
		\|A_\ep^0-\nabla\phi_\ep^0\|_{L^2(\Omega,\R^3)}\leq b^{-1}\|\curl\Bfield\|_{L^2(\Omega,\R^3)}\leq C.
		$$
		Moreover, the hypotheses of Theorem~\ref{thm:main} imply $h_\ex=O(\lep)$. Thus, using \eqref{eq:F-small-final} and $|u|\leq1$, we obtain
		\begin{align*}
			\|\nabla_{\ACouGauge}\uCouGauge\|_{L^2(\Omega,\C^3)}
			&=\|\nabla_Au-ih_\ex(A_\ep^0-\nabla\phi_\ep^0)u\|_{L^2(\Omega,\C^3)}\\
			&\leq b^{-1/2}\sqrt{2\fenweightthreed(u,A)}+Ch_\ex\leq C\lep,
		\end{align*}
		and consequently
		$$
		\|\ACouGauge\|_{W^{1,6}(\Omega,\R^3)}\leq C\lep.
		$$
		Let $\varphi$ be given by Lemma~\ref{lem:fase-Ws} with $s=6$ and $A=\ACouGauge$, and set
		$$
		\AUniGauge\colonequals\ACouGauge-\nabla\varphi,\qquad v\colonequals\uCouGauge e^{-i\varphi}.
		$$
		Then $\AUniGauge$ is in the weighted Coulomb gauge, $(v,\AUniGauge)$ satisfies \eqref{eq:post_splitting_eq}, and $|v|\leq1$. By \eqref{regularitypmin}, \eqref{eq:fase-cota}, and $m(6)=7$,
		$$
		\|\varphi\|_{W^{2,6}(\Omega)}\leq C\big(1+\lep^{M_\infty}\big)^7\lep\leq C\lep^{7M_\infty+1}.
		$$
		Hence
		$$
		\|\AUniGauge\|_{W^{1,6}(\Omega,\R^3)}\leq C\lep^{7M_\infty+1},
		$$
		and Morrey's inequality, together with $k\geq7M_\infty+1$, gives
		\begin{equation}\label{eq:trace-hyp-verified}
			\|\nabla\pmin\|_{L^\infty(\Omega,\R^3)}+\|\AUniGauge\|_{L^\infty(\Omega,\R^3)}+\|\AUniGauge\|_{W^{1,6}(\Omega,\R^3)}\leq C\lep^{7M_\infty+1}\leq C\lep^k.
		\end{equation}
		Thus, the hypotheses fixed before \eqref{eq:definitions_monotonicity} hold for $(v,\AUniGauge)$ with $M_A=k\geq M_\infty$.
		
		\item\label{step:gradient-bound} \textbf{The gradient bound.} In the Coulomb gauge, Proposition~\ref{prop:regularity} gives
		$$
		\|\nabla(\pmin\uCouGauge)\|_{L^\infty(\Omega,\C^3)}\leq C\left(\ep^{-1}+\|\ACouGauge\|_{L^\infty(\Omega,\R^3)}\right)\leq C\ep^{-1}.
		$$
		Since
		$$
		\nabla v=e^{-i\varphi}\big(\nabla\uCouGauge-i\nabla\varphi\,\uCouGauge\big),\qquad \nabla\uCouGauge=\pmin^{-1}\nabla(\pmin\uCouGauge)-\pmin^{-1}\uCouGauge\nabla\pmin,
		$$
		and $|\uCouGauge|\leq1$, we obtain
		\begin{align*}
			\|\nabla v\|_{L^\infty(\Omega,\C^3)}
			&\leq C\left(\|\nabla(\pmin\uCouGauge)\|_{L^\infty(\Omega,\C^3)}+\|\nabla\pmin\|_{L^\infty(\Omega,\R^3)}+\|\nabla\varphi\|_{L^\infty(\Omega,\R^3)}\right)\\
			&\leq C\left(\ep^{-1}+\lep^{M_\infty}+\lep^{7M_\infty+1}\right)\leq\frac{C_1}{\ep}
		\end{align*}
		for $\ep<\ep_0$. Thus, \eqref{eq:grad-hyp} holds.
		
		\item\textbf{The fluctuation in the weighted Coulomb gauge.} As explained in Section~\ref{sec:prelims}, the fluctuation $A$ can be transformed into the weighted Coulomb gauge. We denote the resulting vector potential by
		$$
		\AfluctWeight\colonequals A-\nabla\varphi_A.
		$$
		Similarly, $A_\ep^0-\nabla\phi_\ep^0$ is in the weighted Coulomb gauge and is gauge equivalent to $A_\ep^0$. Since the weighted Coulomb gauge is a vector space, the field
		$$
		\AUniGauge-\AfluctWeight-h_\ex(A_\ep^0-\nabla\phi_\ep^0)
		$$
		also belongs to it. Moreover,
		\begin{align*}
			\AUniGauge-\AfluctWeight-h_\ex(A_\ep^0-\nabla\phi_\ep^0)
			&=\ACouGauge-\nabla\varphi-(A-\nabla\varphi_A)-h_\ex(A_\ep^0-\nabla\phi_\ep^0)\\
			&=-\nabla\varphiCouGauge-\nabla\varphi+\nabla\varphi_A+h_\ex\nabla\phi_\ep^0.
		\end{align*}
		Thus, it is a gradient $\nabla w$ satisfying
		$$
		\diver(\pmin^2\nabla w)=0\quad\text{in }\Omega,\qquad \partial_\nu w=0\quad\text{on }\partial\Omega.
		$$
		Testing this equation with $w$ gives $\nabla w=0$, so $w$ is constant. Therefore,
		$$
		\AUniGauge=\AfluctWeight+h_\ex(A_\ep^0-\nabla\phi_\ep^0),\qquad v=u e^{-i\varphi_A}
		$$
		up to a constant phase. Consequently, $|v|=|u|$ and
		$$
		\nabla v=e^{-i\varphi_A}\big(\nabla_Au+i\AfluctWeight u\big),
		$$
		which implies
		\begin{equation}\label{eq:gradv-fl}
			\frac12|\nabla v|^2\leq|\nabla_Au|^2+|\AfluctWeight|^2\leq b^{-1}\pmin^2|\nabla_Au|^2+|\AfluctWeight|^2.
		\end{equation}
		Let $\Afluctglobfree$ be the representative of the fluctuation in the global divergence-free gauge, so that
		$$
		\curl\Afluctglobfree=\curl A=\curl\AfluctWeight.
		$$
		By Lemma~\ref{lemma:decompositionR3}, there exists $\psi\in H^1(\Omega)$ such that
		$$
		\AfluctWeight=\Afluctglobfree-\nabla\psi.
		$$
		Using $\diver(\pmin^2\AfluctWeight)=0$ and $\AfluctWeight\cdot\nu=0$, we have
		\begin{align*}
			\int_\Omega\pmin^2|\Afluctglobfree|^2
			&=\int_\Omega\pmin^2|\AfluctWeight|^2+2\int_\Omega\pmin^2\AfluctWeight\cdot\nabla\psi+\int_\Omega\pmin^2|\nabla\psi|^2\\
			&=\int_\Omega\pmin^2|\AfluctWeight|^2+\int_\Omega\pmin^2|\nabla\psi|^2\\
			&\geq\int_\Omega\pmin^2|\AfluctWeight|^2.
		\end{align*}
		Therefore, using H\"older's inequality together with \eqref{space1} and \eqref{space2},
		\begin{equation}\label{eq:A-fl-L2}
			\|\AfluctWeight\|_{L^2(\Omega,\R^3)}^2\leq\frac1b\|\Afluctglobfree\|_{L^2(\Omega,\R^3)}^2\leq C\|\Afluctglobfree\|_{L^6(\R^3,\R^3)}^2\leq C\|\curl A\|_{L^2(\R^3,\R^3)}^2\leq C\fenweightthreed(u,A).
		\end{equation}
		
		\item\textbf{Control of $E$ by the weighted free energy.} Here, $e_\ep$ and $b_\ep$ are associated with $(v,\AUniGauge)$, so
		$$
		b_\ep=1-\frac{\ep^2|\AUniGauge|^2}{\pmin^2}.
		$$
		Using $|v|=|u|$, we have
		$$
		b_\ep-|v|^2=(1-|u|^2)-\frac{\ep^2|\AUniGauge|^2}{\pmin^2}.
		$$
		Thus, by \eqref{eq:trace-hyp-verified} and $\pmin^2\geq b$,
		\begin{align*}
			\frac{\pmin^2}{4\ep^2}(b_\ep-|v|^2)^2
			&\leq\frac1b\frac{\pmin^4}{2\ep^2}(1-|u|^2)^2+\frac{\ep^2|\AUniGauge|^4}{2\pmin^2}\\
			&\leq\frac1b\frac{\pmin^4}{2\ep^2}(1-|u|^2)^2+C\ep^2\lep^{4k}.
		\end{align*}
		Combining this estimate with \eqref{eq:gradv-fl}, integrating over $\Omega$, and using \eqref{eq:A-fl-L2}, we obtain
		\begin{equation}\label{eq:E-vs-F}
			\sup_{r>0}E(r)\leq\int_\Omega e_\ep\leq C\fenweightthreed(u,A)+C\ep^2\lep^{4k}.
		\end{equation}
		\item\label{step:clearing-conclusion} \textbf{Conclusion.} By Steps~\ref{step:total-weighted-gauge} and \ref{step:gradient-bound}, the hypotheses of Lemmas~\ref{lem:eps-reg} and \ref{lem:mono-borde} hold for $(v,\AUniGauge)$ with $M_A=k$. Take $\ep_0$ sufficiently small so that
		$$
		\ep<\lep^{-k}<r_1.
		$$
		For $x_0\in\overline\Omega$, apply Lemma~\ref{lem:mono-borde} with $s=\ep$ and $r=\lep^{-k}$. Since $M_A=k$, we have $r\lep^{M_A}=1$, so the exponential factor in \eqref{eq:mono} is uniformly bounded. Moreover,
		$$
		\Phi(\lep^{-k})=\lep^kE(\lep^{-k}),\qquad K_\ep=\ep^2\lep^{5k}.
		$$
		Recalling that $|u|=|v|$, we obtain
		\begin{align*}
			(1-|u(x_0)|)^5
			&=(1-|v(x_0)|)^5\\
			&\overset{\eqref{eq:eps-reg}}{\leq}C\big(\Phi(\ep)+\ep^4\lep^{4k}\big)\\
			&\overset{\eqref{eq:mono}}{\leq}C\big(\lep^kE(\lep^{-k})+K_\ep+\ep^4\lep^{4k}\big)\\
			&\overset{\eqref{eq:E-vs-F}}{\leq}C\big(\lep^k\fenweightthreed(u,A)+\ep^2\lep^{5k}+K_\ep+\ep^4\lep^{4k}\big)\\
			&\leq C\big(\lep^k\fenweightthreed(u,A)+\ep^2\lep^{5k}\big)\\
			&\overset{\eqref{eq:F-small-final}}{=}o(\lep^{-1}).
		\end{align*}
		This chain of estimates implements the clearing-out strategy: \eqref{eq:eps-reg} controls the pointwise defect by the normalized energy at scale $\ep$, \eqref{eq:mono} transports energy control from the larger scale $\lep^{-k}$ down to scale $\ep$, \eqref{eq:E-vs-F} relates the resulting local energy to the free energy of the fluctuation, and \eqref{eq:F-small-final} provides its required smallness. Since the estimate is uniform in $x_0$, it follows that
		$$
		\|1-|u|\|_{L^\infty(\Omega)}=o(1).
		$$
		Finally, since $|\u|=\pmin|u|$ and $\|\pmin\|_{L^\infty(\Omega)}\leq1$,
		$$
		\|\pmin-|\u|\|_{L^\infty(\Omega)}\leq\|\pmin\|_{L^\infty(\Omega)}\|1-|u|\|_{L^\infty(\Omega)}\leq\|1-|u|\|_{L^\infty(\Omega)}=o(1).
		$$
		This proves \eqref{eq:clearing-final}.
		
		\item \textbf{Improved vorticity and free-energy estimates.} By Step~\ref{step:clearing-conclusion}, $|u|\geq1/2$ in $\Omega$ for $\ep$ sufficiently small. We may therefore apply \cite[Proposition A.1]{Rom-CMP}, which gives
		\begin{equation}\label{eq:improved-vorticity-lip}
			\|\mu(u,A)\|_{(C_T^{0,1}(\Omega,\R^3))^*}\leq C\ep \fenthreedOmega(u,A),
		\end{equation}
		where
		$$
		\fenthreedOmega(u,A)\colonequals\frac12\int_\Omega\left(|\nabla_Au|^2+\frac1{2\ep^2}(1-|u|^2)^2+|\curl A|^2\right).
		$$
		Since $b\leq\pmin^2\leq1$, we have
		\begin{equation}\label{eq:unweighted-vs-weighted}
			\fenthreedOmega(u,A)\leq b^{-2}\fenweightthreed(u,A).
		\end{equation}
		In addition, from \cite[Lemma 8.1]{roman-vortex-construction} we have
		\begin{equation}\label{vorticityestimatecontinuous}
			\|\mu(u,A)\|_{(C^0(\Omega,\R^3))^*}\leq C\fenthreedOmega(u,A).
		\end{equation}
		Interpolating between \eqref{eq:improved-vorticity-lip} and \eqref{vorticityestimatecontinuous}, as explained in \cite[Section 8]{roman-vortex-construction}, we obtain, for every $\gamma\in(0,1]$,
		\begin{equation}\label{eq:improved-vorticity-holder}
			\|\mu(u,A)\|_{(C_T^{0,\gamma}(\Omega,\R^3))^*}\leq C_\gamma\ep^\gamma \fenthreedOmega(u,A).
		\end{equation}
		Fix $\gamma\in(0,1)$. By \eqref{estimateforB0holder}, \eqref{eq:improved-vorticity-holder}, and \eqref{eq:unweighted-vs-weighted},
		$$
		h_\ex\left|\int_\Omega\mu(u,A)\cdot\Bfield\right|\leq Ch_\ex\ep^\gamma\lep^{\McurlBfield}\fenthreedOmega(u,A)\leq Ch_\ex\ep^\gamma\lep^{\McurlBfield}\fenweightthreed(u,A).
		$$
		Since $h_\ex=O(\lep)$, the coefficient on the right-hand side tends to zero and can be absorbed into the left-hand side of \eqref{eq:upper-split}. Using also \eqref{eq:r-bound}, we obtain
		$$
		\fenweightthreed(u,A)\leq Ch_\ex^2\ep\lep^{2\McurlBfield}\fenweightthreed(u,A)^{1/2}.
		$$
		Consequently,
		\begin{equation}\label{estfen}
			\fenweightthreed(u,A)\leq Ch_\ex^4\ep^2\lep^{4\McurlBfield}\leq C\ep^2\lep^{4+4\McurlBfield}.
		\end{equation}
		In particular, by \eqref{eq:unweighted-vs-weighted},
		\begin{equation}\label{eq:improved-unweighted-free-energy}
			\fenthreedOmega(u,A)\leq C\ep^2\lep^{4+4\McurlBfield}=o(\ep).
		\end{equation}
		Combining \eqref{eq:improved-vorticity-holder} and \eqref{eq:improved-unweighted-free-energy}, we also obtain, for every $\gamma\in(0,1]$,
		\begin{equation}\label{estvort}
			\|\mu(u,A)\|_{(C_T^{0,\gamma}(\Omega,\R^3))^*}\leq C_\gamma\ep^{2+\gamma}\lep^{4+4\McurlBfield}.
		\end{equation}
		Finally, from \eqref{estfen} and \eqref{estvort}, we obtain \eqref{eq:improved-conclusions}.		
	\end{enumerate}
	The theorem is thus proved.
\end{proof}

\begin{remark}[Connection with the metastable Meissner branch]
	The estimate \eqref{eq:improved-unweighted-free-energy} gives
	$$
	\fenthreedOmega(u,A)=o(\ep),
	$$
	which places the global minimizer in the uniqueness regime considered in \cite{rom2026metastab}. Therefore, under the additional assumptions of that work, the Meissner-type local minimizer constructed in \cite{rom2026metastab} is in fact the unique global minimizer, up to gauge transformations, below the first critical field.
\end{remark}

\bibliography{references3DPinning}
\end{document}